\documentclass{article}[11]

\usepackage{graphics,latexsym}
\usepackage{graphicx}
\usepackage{amsmath, amssymb,amsthm}

\usepackage{color}

\begin{document}

\title{Optimal Collocation Nodes for Fractional Derivative Operators}

\author{Lorella Fatone$^{1}$, Daniele Funaro$^{2}$}

\maketitle

\centerline{$^1$\small Dipartimento di Matematica  e Informatica }
\centerline{\small Universit\`a di Camerino, Via Madonna delle Carceri 9, 62032
Camerino (Italy)} 
\medskip

\centerline{$^2$\small Dipartimento di Fisica, Informatica e Matematica }
\centerline{\small Universit\`a di Modena e Reggio Emilia, Via Campi 213/B, 41125
Modena (Italy)} 
\medskip

%\centerline{\bf \today}

\medskip

\begin{abstract} {Spectral discretizations of fractional derivative
operators are examined, where the approximation basis is related to the set of
Jacobi polynomials. The pseudo-spectral method is implemented by assuming
that the grid, used to represent the function to be differentiated,
may not be coincident with the collocation grid. The new option opens the way
to the analysis of alternative techniques and the search of optimal
distributions of collocation nodes, based on the operator to be
approximated. Once the initial representation grid has been chosen,
indications on how to recover the collocation grid are provided,
with the aim of enlarging the dimension of the approximation space.
As a results of this process, performances are
improved. Applications to  fractional type advection-diffusion equations, 
and comparisons in terms of accuracy and efficiency are made.
As shown in the analysis, special choices of the nodes can also suggest
tricks to speed up computations. 
}
\end{abstract}

\vspace{.2cm}
\par\smallskip

\vspace{.2cm}
\par\smallskip

\maketitle

%%%%%%%%%%%%%%%%%%%%%%%%%%%%%%%%%%%%%%%%%%%
\section{Aim of the paper}\label{sec1}
%%%%%%%%%%%%%%%%%%%%%%%%%%%%%%%%%%%%%%%%%%%
Boundary-value problems involving derivatives of fractional
order have found increasing interest in the last years.
They emerge in a large number of applications, ranging from quantum 
mechanics to mechanical engineering, chemistry or economics.
The literature offers a wide collection of papers.
Some references in alphabetical order are for instance: 
\cite{barkai}, \cite{benson},  \cite{chaves}, \cite{ervin1}, \cite{ervin2},
\cite{henry}, \cite{magin}, \cite{maina}, \cite{metz1},
\cite{metz2}, \cite{pablo}, \cite{podlu}, \cite{sugi1}, \cite{sugi2}, \cite{zas}. 
Analytical solutions via Laplace, Fourier, or Mellin transforms have
been proposed in several of the above mentioned papers.
\par\smallskip

In the framework of numerical approximations, investigations have developed 
along different paths, including finite-difference methods (\cite{deng2}, \cite{dieth1}, 
\cite{dieth2}, \cite{goren}, \cite{liu2}, \cite{lubi}, \cite{meer1}, \cite{meer2},
\cite{metz2}, \cite{sugi1}, \cite{sun})
and finite element methods (\cite{deng1}, \cite{ervin1}, \cite{ervin2}, \cite{zhang}).
More recently, high-order techniques have also been employed. These 
involve the use of spectral Galerkin methods (\cite{lixu}, \cite{lixu2} \cite{sugi1},
\cite{lin}) or spectral collocation (\cite{khader2}, \cite{tian}, \cite{zakar1}, \cite{zakar2}).
The present paper deals with the last subject.
\par\smallskip

Using as approximation basis the set of Jacobi polynomials, pseudo-spectral
discretizations of fractional derivative operators are introduced and examined.
The idea is to ameliorate the methods recently proposed in \cite{tian} and \cite{zakar2}.
To this end, a suitable techniques is suggested, where
the grid used to represent the function to be differentiated,
is not necessarily coincident with the collocation grid. This option
was studied in \cite{miolibro2}, \cite{fatfun} and \cite{funsup} in the framework of standard 
partial differential equation and in \cite{funint} for integral type equations.
The scope of using two grids is to enlarge, through a procedure named
{\sl superconsistency}, the dimension of the approximation space.
The result is an improvement of the overall performances of the 
method, with very little additional cost. 
\par\smallskip

Asymptotically, i.e. when the
number of nodes increases, having different sets for the representation
and collocation nodes does not bring to drastic differences. Nevertheless, for lower 
degree approximating polynomials the gain may be very impressive.
The methodology is then appropriate for stiff problems, where the 
number of degrees of freedom used for the discretization is
still low in order to achieve high accuracy. Examples of this kind
are transport-diffusion equations with dominating advective terms.
If the polynomial degree is too small to resolve boundary layers,
the approximate solution may be very rough. According to 
\cite{miolibro2} and \cite{fatfun}, the adoption of a suitable collocation
grid brings to excellent improvements. For this reason, in the last
section of this paper, we examine a transport-diffusion equation, where
the operator contains fractional derivatives. We compare
different collocation procedures showing that our approach is
actually competitive.
\par\smallskip

Our discussion also involves a review of the construction of the
approximation matrices. In some special cases, we will be able to come
out with an explicit expression of the entries of the linear
discrete operators. Usually, these quantities are instead computed by
introducing further approximation.
\par\smallskip

We complete this short introduction with some preliminary definitions. 
We are concerned with computing fractional derivatives in the interval 
$\left[ -1,1 \right]$. To this end we work with the Riemann-Liouville fractional 
operator of order $\sigma$:
\begin{equation}\label{deri}
(D^\sigma f)(x)=\frac{1}{\Gamma (1-\sigma)}~\frac{d}{dx}\int_{-1}^x
\frac{f(s) ds}{(x-s)^\sigma}, \quad \quad x>-1.
\end{equation}
Here $0 <\sigma <1$  is the derivative order and $\Gamma $ denotes the Euler gamma function.
 Other versions of 
fractional derivative operators, such as the Caputo's, are
available. They are all connected by simple relations, so that,
what we are going to develop in the coming sections can be
easily extended to other cases.
For a general survey of fractional calculus see, for example, \cite{podlu}.
\par\smallskip

As we said above, our interest is mainly focused on high-order approximation techniques.
We will mainly use collocation type methods based on the zeros of Jacobi polynomials. 
For the reader's convenience, we briefly review some basic and remarkable properties of Jacobi polynomials
First of all, we recall that Jacobi polynomials are denoted by
$P_n^{\alpha ,\beta}$ where $n\geq 0$ is the degree and $\alpha >-1$
and $\beta >-1$ are given parameters. 
For $n\geq 0$, the $n$-th Jacobi polynomial satisfies the following Sturm-Liouville 
eigenvalue problem in $\left[ -1,1 \right]$:
\begin{equation}\label{slp}
(1-x^2)\frac{d^2 P_n^{\alpha ,\beta}}{dx^2}-\big( (\alpha +\beta +2)x
+\alpha -\beta \big)\frac{d P_n^{\alpha ,\beta}}{dx} +n(n+\alpha +\beta +1)
P_n^{\alpha ,\beta}=0.
\end{equation}
\smallskip
Jacobi polynomials are characterized
by the orthogonality relation:
\begin{equation}\label{orto}
\int_{-1}^1 P_n^{\alpha ,\beta}(x)P_k^{\alpha ,\beta}(x)(1-x)^\alpha
(1+x)^\beta dx=0, \quad \quad {\rm if}\, k\not =n.
\end{equation}
Moreover, one has for $n\geq 1$:
$$
\int_{-1}^1 \left[ P_n^{\alpha ,\beta}(x)\right]^2(1-x)^\alpha
(1+x)^\beta dx
$$
\begin{equation}\label{orton}
=\frac{2^{\alpha +\beta +1}}{n!~(2n+\alpha +\beta +1)}
\frac{\Gamma (n+\alpha +1)\Gamma (n+\beta +1)}{\Gamma (n+\alpha +\beta+1)}, 
\end{equation}
\smallskip
For $n\geq 1$, a very useful relation is:
\begin{equation}\label{derip}
\frac{d}{dx}\left[ P_n^{\alpha, \beta}\right]=\frac{n+\alpha +\beta +1}{2}
P_{n-1}^{\alpha +1, \beta +1}.
\end{equation}
We finally recall that, starting from: 
\begin{equation}\label{recursive1}
P_0^{\alpha, \beta}(x)=1, \quad \qquad P_1^{\alpha, \beta}(x)=\frac12 \left(  \alpha+\beta+2 \right) x +\frac12 (\alpha-\beta),
\end{equation}
\smallskip
higher degree Jacobi polynomials can be determined using the 
following recurrence  relation:
\begin{equation}\label{recursive2}
P_n^{\alpha, \beta} (x)=\left(  a_{n} x +b_{n}  \right) P_{n-1}^{\alpha, \beta} (x)+c_{n} P_{n-2}^{\alpha, \beta} (x), \qquad \forall  n\ge 2 ,
\end{equation}
\smallskip
where:
\begin{eqnarray}
&& a_{n} = \frac{ \left(2n+  \alpha+\beta \right)  \left(  2n+ \alpha+\beta-1 \right)     }{2n  \left(  n+ \alpha+\beta \right)}    ,    \label{reca}\\[3mm]
&& b_{n} = \frac{ \left(\alpha^{2}-\beta^{2} \right)  \left(  2n+ \alpha+\beta-1 \right)     }{2n  \left(  n+ \alpha+\beta \right) \left( 2 n+ \alpha+\beta  -2\right) }    ,    \label{recb}\\[3mm]
&&c_{n} =  - \frac{ \left(  n+ \alpha-1 \right)  \left(  n+\beta-1 \right)     \left(  2n+ \alpha+\beta \right)     }{n  \left(  n+ \alpha+\beta \right) \left( 2 n+ \alpha+\beta  -2\right) }    , \qquad \forall  n\ge 2.\label{recc}
\end{eqnarray}
Ultraspherical polynomials are Jacobi polynomials where $\alpha=\beta$. 
Legendre polynomials are ultraspherical polynomials with $\alpha=\beta=0$.
In order to simplify the notation, we set   $P_n=P_n^{0 ,0}$.
Chebyshev  polynomials  (of the first kind) are related to the 
ultraspherical polynomials with $\alpha=\beta=-\frac12$. In fact, they are defined by:
\begin{equation}\label{Cheb}
T_{n}(x)=\frac{\left(   n! ~ 2^{n} \right)^{2}}{\left( 2n \right)!} \, 
P_n^{-\frac12 ,-\frac12 }(x), \quad \quad n\ge 0.
\end{equation}
For a complete survey  of the properties of Jacobi, Legendre and Chebyshev polynomials,   
as well as other commonly used families of orthogonal polynomials,
we refer for instance to \cite{Szego} and \cite{miolibro}.
\par\smallskip

The most important relation linking Jacobi polynomials with fractional
derivatives is represented by the following  equation (see \cite{Askey}):
$$
\frac{\Gamma (\beta+ \mu +1) P_n^{\alpha -\mu,\beta+\mu}(-1)}{\Gamma (\beta +1)\Gamma (\mu )
P_n^{\alpha ,\beta}(-1)}\int_{-1}^x\frac{(1+s)^\beta P_n^{\alpha ,\beta}(s)}
{(x-s)^{1-\mu}} \, ds=(1+x)^{\beta +\mu} P_n^{\alpha -\mu,\beta+\mu}(x),
$$
\begin{equation}\label{auto}
\hfill \, \, \,\, 0<\mu <1, \,\,\, x \in \left[ -1,1 \right].
\end{equation}
\par\smallskip
An interesting version of (\ref{auto}) is obtained  when $\alpha =\mu$, $\beta =-\mu$.
With this choice one has:
$$
\frac{P_n(-1)}{\Gamma (1-\mu)\Gamma (\mu )
P_n^{\mu ,-\mu}(-1)}\int_{-1}^x\frac{(1+s)^{-\mu} P_n^{\mu ,-\mu}(s)}
{(x-s)^{1-\mu}}\,ds= P_n(x), 
$$
\begin{equation}\label{autom}
\hfill  0<\mu <1, \,\,\, x \in \left[ -1,1 \right].
\end{equation}
The analysis  carried out in this paper is mainly based on the above formula, but it is clear that
straightforward generalizations are possible by using the
full potentiality of (\ref{auto}).
We are now ready to study approximations of the operator $D^\sigma$, $0 <\sigma <1$.
\par\smallskip

%%%%%%%%%%%%%%%%%%%%%%%%%%%%%%%%%%%%%%%%%%%
\section{A collocation method for ${\bf D}^{\boldsymbol\sigma} $, ${\bf 0 <}
{\boldsymbol\sigma}{\bf <1}$}\label{sec2}
%%%%%%%%%%%%%%%%%%%%%%%%%%%%%%%%%%%%%%%%%%%

We assume that the nodes $x_j\in \left[ -1,1\right]$, $j=0,1,\ldots,N$, are given for some
integer $N\geq 2$.  Their explicit expression will be examined later.
From now on we set: $x_0=-1$. Afterwards, we suppose to have a function $u_N$,  
satisfying $u_N(-1)=0$ and depending on $N$ degrees of freedom.
Similarly to what has been done in  \cite{zakar2},  given $0<\mu <1$, we suppose that
$u_N$ is expanded in the following Lagrange basis:
\begin{equation}\label{svil}
u_N(x)=\sum_{j=1}^N u_N(x_j) H_j(x), \quad \quad 0<\mu <1, \,\,\, x \in \left[ -1,1 \right],
\end{equation}
where  the basis elements $H_j$, $j=1,2,\ldots,N$, are defined as follows:
\begin{eqnarray}
&&H_j(x)=\left(\frac{x_j+1}{x+1}\right)^\mu \prod_{{k=0}\atop{k\not =j}}^N
\left(\frac{x-x_k}{x_j-x_k}\right)
=
\left(\frac{x+1}{x_j+1}\right)^{1-\mu} \prod_{{k=1}\atop{k\not =j}}^N
\left(\frac{x-x_k}{x_j-x_k}\right),\nonumber\\[3mm]
&&\hskip6truecm \quad0<\mu <1, \,\,\, x \in \left[ -1,1 \right].\label{acca}
\end{eqnarray}
Indeed, we have the Kronecker delta property $H_j(x_m)=\delta_{jm}$, 
$j,m=1,2,\ldots,N$.
Basically, we have to deal with a polynomial of degree $N$, suitably
corrected at the point $x_{0}=-1$. Despite this initial setting,
our approach is going to be different from the one followed
in  \cite{zakar2}.
\par\smallskip

By the linearity of the fractional differential operator $D^\sigma$,  $0 <\sigma <1$, we are allowed to
write:
\begin{equation}\label{dsvil}
(D^\sigma u_N)(x)=\sum_{j=1}^N u_N(x_j) (D^\sigma H_j)(x), \quad \quad 0<\sigma <1, 
\,  x \in \left[ -1,1 \right],
\end{equation}
so that we need to evaluate the effect of applying $D^\sigma$, $0<\sigma <1$, 
to each element of the basis $H_j$, $j=1,2,\ldots,N$. 
\par\smallskip

In order to use formula  (\ref{autom}),  we are going to represent the elements $H_j$, $j=1,2,\ldots,N$,
in the following equivalent form:
\begin{equation}\label{svilh}
H_j(x)=(x+1)^{-\mu}\sum_{n=1}^N c(j,n) \left[ P_n^{\mu ,-\mu}(x)
-P_n^{\mu ,-\mu}(-1) \right].
\end{equation}

Note that the right-hand side  of (\ref{svilh}) actually vanishes for $x_{0}=-1$,
because $H_j(x)\approx (x+1)^{1-\mu}$ in the neighborhood of that point.
\par\smallskip

The difficulty is to pass from one representation of $u_N$
to the other, that is to pass from  the representation (\ref{svil})  of $u_N$ with  $H_j$, $j=1,2,\ldots,N$,
defined in  (\ref{acca})  to the representation (\ref{svil}) of $u_N$  with  $H_j$, $j=1,2,\ldots,N$,
defined in (\ref{svilh}). 
\par\smallskip

Searching for explicit formulas is rather
cumbersome and may lead to an expensive and  ill-conditioned
algorithm, especially if one goes through to the monomial
basis $x^k$, $k=1,2,\ldots,N$ (which is the procedure followed in \cite{tian}). 
The solution we suggest is to solve a simple linear
system. In fact, by observing that $H_j(x_m)=\delta_{jm}$, $j,m=1,2,\ldots,N$, 
one can evaluate equation (\ref{svilh}) at $x=x_m$, $m=1,2,\ldots,N,$
obtaining:
\begin{equation}\label{svilhm}
(x_m+1)^{-\mu}\sum_{n=1}^N c(j,n)  \left[ P_n^{\mu ,-\mu}(x_m)
-P_n^{\mu ,-\mu}(-1)\right]= \delta_{jm}. 
\end{equation}
By introducing the $N\times N$ matrix:
\begin{equation}\label{mat}
A_N=\{ a_{j,n}\}=\{ (x_j+1)^{-\mu}\left[ P_n^{\mu ,-\mu}(x_j)
-P_n^{\mu ,-\mu}(-1)\right]\},   
\end{equation}
where $j,n=1,2,\ldots,N$,
we get the expansion coefficients $c(j,n)$  as the entries of the matrix  $A_N^{-1}$.
This computation is relatively cheap and
the condition number of $A_N$ is quite acceptable, as one can
check, for example,  by examining Table \ref{table1}. This test is related to the distribution of nodes:
\begin{equation}\label{nodiCheb}
x_j=-\cos (j\pi /N), \quad j=1,2,\ldots,N.
\end{equation}
 The nodes $x_j$, $j=1,2,\ldots,N$  defined in (\ref{nodiCheb}) are the zeros of $T^\prime_N$
(the derivative  of the $N$-th Chebyshev 
polynomial)    with the addition of the point $x_N=1$.
\par\smallskip

 Table \ref{table1}  shows that the growth of the condition number of $A_N$ 
is clearly proportional to $N$ in this case. In some circumstances the 
evaluation of $c(j,n)$, $j,n=1,2,\ldots,N$,  is straightforward, as it will be checked in
Section 5.
\par\smallskip

\begin{table}
\centerline{
{\begin{tabular}{@{}ccclccc}\hline
 &&$N$ &&  cond$(A_{N})$&& \\
\hline
 &&  5   && 3.7240 &&\\
 &&  10 && 8.9481 &&\\
 &&  20 && 19.0645 &&\\
 &&  50 &&  50.3533 &&\\
 &&  100 &&  103.4209 &&\\
 \hline
  \end{tabular}}
  }\vspace{-0.5cm}
  \bigskip\caption{\small Condition number of $A_N$ 
when $x_j=-\cos (j\pi /N)$, $j=1,2,\ldots,N$ and $\mu=0.5$.}
\label{table1}
\end{table}

We are now ready to compute the fractional derivative of all $H_j$,
$j=1,2,\ldots,N$. For this, we require that $\mu =1-\sigma$. Therefore using (\ref{autom}) we obtain:
$$
(D^\sigma H_j)(x)=(D^{1-\mu} H_j)(x)=\frac{1}{\Gamma (\mu)}~\frac{d}{dx}\int_{-1}^x
\frac{H_j(s) ds}{(x-s)^{1-\mu}}
$$
$$
=\frac{1}{\Gamma (\mu)}
\sum_{n=1}^N c(j,n)\frac{d}{dx}\int_{-1}^x\frac{(1+s)^{-\mu} P_n^{\mu ,-\mu}(s)}
{(x-s)^{1-\mu}}
$$
$$
=\sum_{n=1}^N\frac{\Gamma (1-\mu )P_n^{\mu ,-\mu}(-1)}{P_n(-1)}
~c(j,n) P_n^\prime (x)
$$
\begin{equation}\label{derih}
=\sum_{n=1}^N  \frac{\Gamma (n-\mu +1)}{n!} 
~c(j,n) P_n^\prime (x).
\end{equation}
In the last passage we used the following relation (see \cite{Szego}):
\begin{equation}\label{jac}
 P_n^{\alpha ,\beta}(-1)=(-1)^n\left({{n+\beta }\atop{n}}\right)=(-1)^n~
 \frac{\Gamma (n+\beta +1)}{n! ~\Gamma (\beta +1)},
\end{equation}
where $\alpha=\mu$ and $\beta=-\mu$.
In particular, one has $P_n(-1)=(-1)^n$.
In computing (\ref{derih}) we eliminated a term. In truth, that
was actually zero. Indeed, by recalling that $ P_0^{\alpha ,\beta}(x)=1$,
$\forall x$, we have:
$$
\frac{d}{dx}\int_{-1}^x\frac{(1+s)^{-\mu} P_n^{\mu ,-\mu}(-1)}
{(x-s)^{1-\mu}} ds=P_n^{\mu ,-\mu}(-1)\frac{d}{dx}\int_{-1}^x\frac{(1+s)^{-\mu} 
P_0^{\mu ,-\mu}(s)}{(x-s)^{1-\mu}}
$$
\begin{equation}\label{ter0}
=P_n^{\mu ,-\mu}(-1)\frac{\Gamma (1-\mu )\Gamma (\mu )P_0^{\mu ,-\mu}(-1)}{P_0(-1)}
P_0^\prime (x)=0.
\end{equation}
Finally, given $0 <\sigma <1$, by previously computing the coefficients $c(j,n)$, we can 
assemble the $N\times N$ fractional derivative matrix $D_N^\sigma$  by setting:
\begin{equation}\label{ddif}
D_N^\sigma =\{ d^\sigma_{i,j}\}=\sum_{n=1}^N  \frac{\Gamma (n-\mu +1 )
}{n!}~ c(j,n) P_n^\prime (z_i), \quad \quad i, j=1,2,\ldots,N. 
\end{equation}
In (\ref{ddif}) the collocation nodes $z_i$, $i=1,2,\ldots,N$ do not
necessarily coincide with the nodes $x_j$, $j=1,2,\ldots,N$.
\par\smallskip

In the end, for a given $0 <\sigma <1$, let us suppose to have the   fractional  differential 
problem:
\begin{equation}\label{pru}
D^\sigma u=g~~~~~~~~~~ u(-1)=0,
\end{equation}
with given right-hand side $g$.
We then propose to approximate (\ref{pru}) by the following discrete
problem:
\begin{equation}\label{prud}
D_N^\sigma u_N=g_N~~~~~~~~~~ u_N(-1)=0,
\end{equation}
where  $D_N^\sigma$ replaces $D^\sigma$,  $ u_N$ is the discrete solution and $g_N$  is the  interpolant of $g$ at the collocation nodes.
\par\smallskip

In the following sections we will specify how to properly choose 
representation and collocation nodes in order to achieve optimal
results.

%%%%%%%%%%%%%%%%%%%%%%%%%%%%%%%%%%%%%%%%%%%
\section{Higher-order fractional derivative operators}\label{sec3}
%%%%%%%%%%%%%%%%%%%%%%%%%%%%%%%%%%%%%%%%%%%

Fractional derivatives of order greater that one can be  obtained by
composition. In particular  fractional derivatives of order $1+ \sigma$,  $0 <\sigma <1$, 
can be computed as follows:

\begin{equation}\label{derp}
D^{1+\sigma}=D D^{\sigma}=D^{\sigma}D  ,
\end{equation}
where $D=D^{1}$ is the standard first derivative operator.
Similarly, one can handle derivatives of the form $D^{k+\sigma}$, $0 <\sigma <1$,
where $k$ is an integer such that $k\geq 1$.
\par\smallskip

Let us note that, given an integer $k\geq 1$, standard derivatives  $D^k $ of order $k$ of the basis functions 
$H_j$ in (\ref{svilh}) are evaluated as follows:
\begin{equation}\label{svilhf}
(D^k H_j)(x)=\sum_{n=1}^N c(j,n)\frac{d^k}{dx^k}\left[(x+1)^{-\mu} (P_n^{\mu ,-\mu}(x)
-P_n^{\mu ,-\mu}(-1))\right], \, k\geq 1.
\end{equation}

For example in Section \ref{sec6}, we will approach a differential equation involving the operator
$-D^2+KD^\sigma $, where $K$ is a  given constant.
By taking $k=2$ in (\ref{svilhf}) and combining (\ref{derih}) and (\ref{svilhf}), 
one  gets for $j=1,2,\ldots,N$:
$$
(-D^2 H_j+KD^\sigma H_j )(x)=\sum_{n=1}^N  c(j,n)\left( 
-\frac{d^2}{dx^2}\left[ (x+1)^{-\mu} P_n^{\mu ,-\mu}(x)
\right]\right.
$$
\begin{equation}\label{dercom2}
\left. +\mu (\mu+1) (x+1)^{-\mu -2}P_n^{\mu ,-\mu}(-1)
+K \frac{\Gamma (n-\mu +1)}{n!} P_n^\prime (x)\right).
\end{equation}

Given an integer $k\geq 1$ more involved fractional differential operators, such as
$-D^{k+\sigma_1}+ K D^{\sigma_2}$ with $0 <\sigma_{1} <1$, $0 <\sigma_{2} <1$ and $\sigma_1\not =\sigma_2$,
can be approached in the way described above starting from the following relations:
\begin{equation}\label{des1}
-D^{k+\sigma_1}+ K D^{\sigma_2}= D^{\sigma_2}( -D^{k+\sigma}+K )~~~~~\quad
\sigma =\sigma_1 -\sigma_2, ~~ \sigma_1\geq \sigma_2 , 
\end{equation}
\begin{equation}\label{des2}
-D^{k+\sigma_1}+ K D^{\sigma_2}= D^{\sigma_1}( -D^k+K D^\sigma)~~~~~\quad
\sigma =\sigma_2 -\sigma_1, ~~ \sigma_2\geq \sigma_1.
\end{equation}
The problem of discretizing these operators becomes however too 
technical and will not be discussed in the present paper.
\par\smallskip

We conclude this section with some explicit examples.
Given an integer $k\geq 1$,  we can build the discretization   
$D^{k+\sigma}_N$ of $D^{k+\sigma}$, $0 <\sigma <1$, by taking for example, as representation
nodes $x_j$, $j=1,2,\ldots,N$, the zeros of the derivative  of the Chebyshev 
polynomial $T_N$, with the additional point $x_N=1$, i.e. the nodes 
defined in  (\ref{nodiCheb}). Moreover, we suppose for
the moment that the collocation nodes  coincide with the representation nodes; 
in other words we assume: $z_i=x_i$, $i=1,2,\ldots,N$.
\par\smallskip

\begin{figure}[h!]
\vspace{.1cm}
\centering
{\includegraphics[width=5.8cm,height=4.8cm]{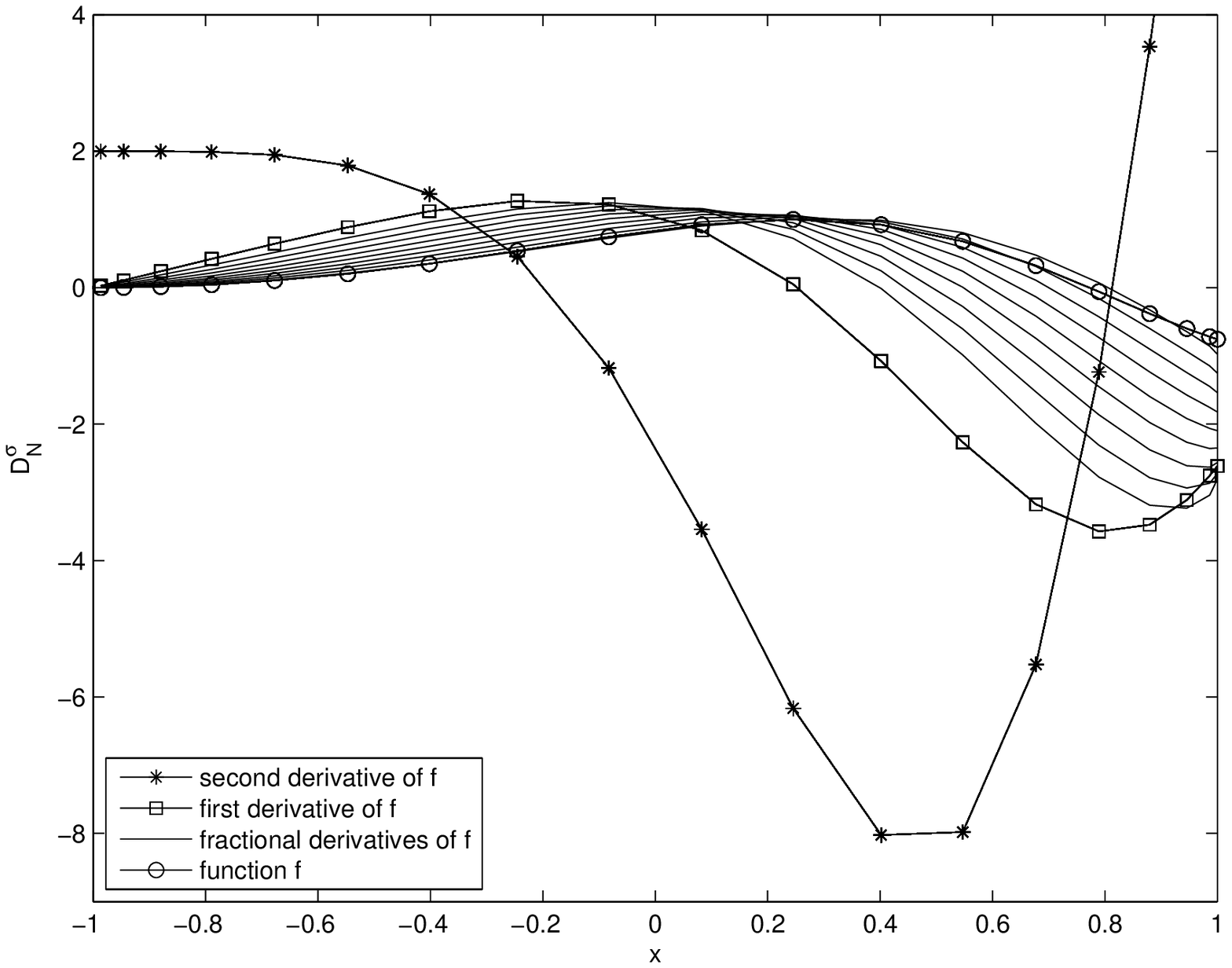}
\hspace{.3cm}\includegraphics[width=5.8cm,height=4.8cm]{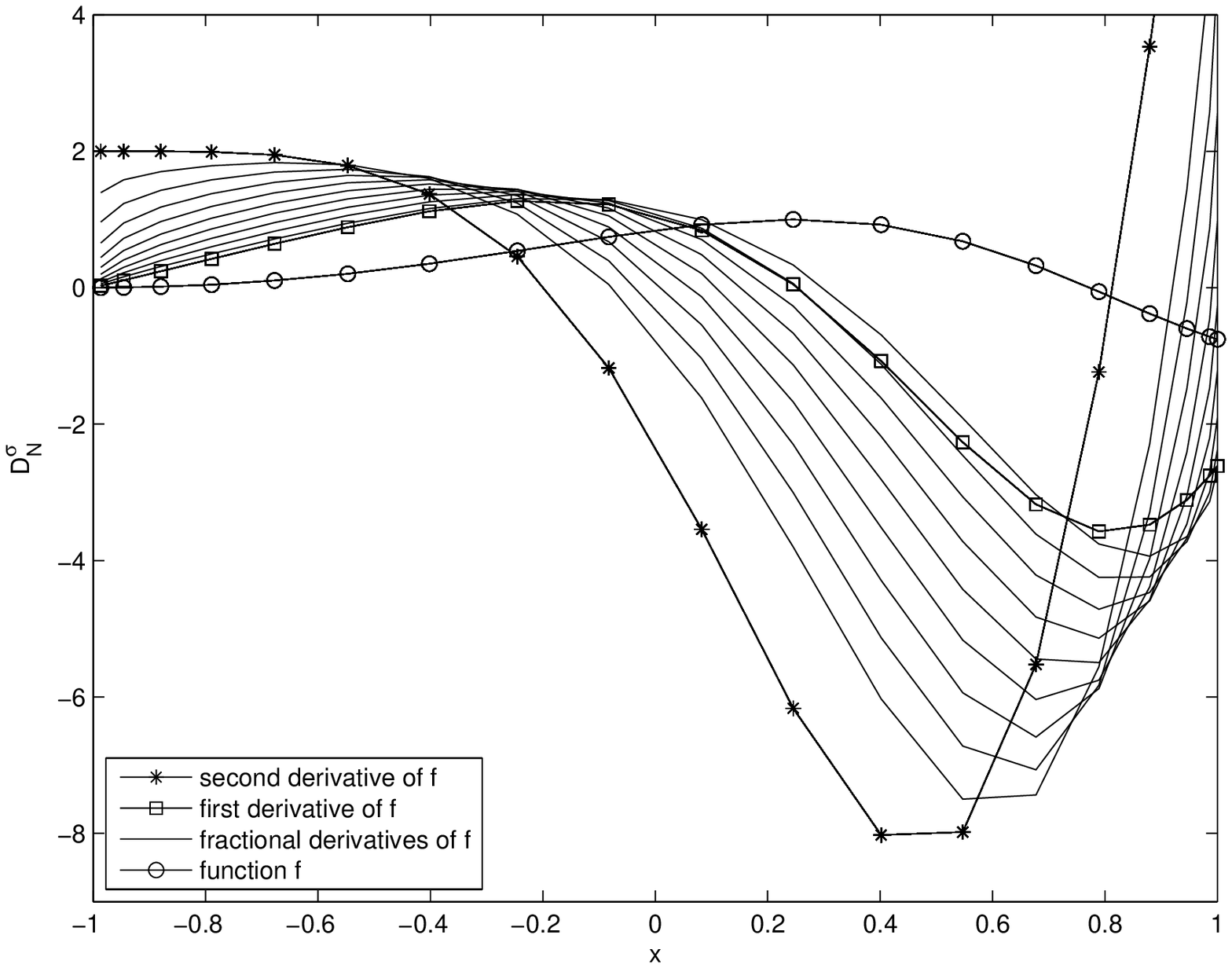}}
\caption{\small \sl Fractional derivative approximations $D^\sigma_N$ of the function
$f(x)=\sin (x+1)^2$ for $N=19$. Here $\sigma$ varies from $0.1$ to $0.9$, step $0.1$
(left), and from $1.1$ to $1.9$, step $0.1$ (right).}  \label{fig1}
\end{figure}

In Figure  \ref{fig1}, we show the results of some tests. Fractional
derivatives of the function $f(x)=\sin (x+1)^2$ are computed for
a given $N$ and various choices of $\sigma$, $0<\sigma <1$. The decay of $f$ near
the point $x_{0}=-1$ is quadratic. This allows for a rather good calculation
of the derivatives up to the order one. The decay of $f^\prime$ is
just linear and this creates a kind of boundary layer near
the point $x=-1$. The reason for this behavior can be attributed to
the decision of representing $u_N$ through the basis in (\ref{acca}),
where the parameter $\mu =1-\sigma$ dictates the decay rate of the discrete
fractional derivative at the point $x=-1$. In order to handle these
specific situations, a less tamed basis should be
constructed on purpose, though this is not a subject we shall deal with.
\par\smallskip

%%%%%%%%%%%%%%%%%%%%%%%%%%%%%%%%%%%%%%%%%%%
\section{Choice of the collocation nodes}\label{sec4}
 %%%%%%%%%%%%%%%%%%%%%%%%%%%%%%%%%%%%%%%%%%%
 
As we mentioned in the previous section, we are not obliged
to choose the set of collocation nodes equal to that used
to represent the solution. This observation suggests a series 
of experiments with different combinations of nodes.
Although the choices can be infinite, we will concentrate
our attention on some meaningful cases.
\par\smallskip

First of all, we note that, thanks to (\ref{derip}), the derivative of the Jacobi polynomial
$P^{\alpha -1 ,\beta -1 }_N$ is proportional to $P^{\alpha ,\beta }_{N-1}$.
As suggested in \cite{zakar2}, a framework providing  very good 
performances is the one where the representation nodes $x_j$, $j=0,2,\ldots,N$,  are the zeros of 
$P^{\alpha ,\beta }_{N-1}$, with the addition of the points 
$x_0=-1$ and $x_N=1$ (see (\ref{nodiCheb}) concerning the Chebyshev case).

\par\smallskip

Systematically, the collocation nodes are chosen such that $x_j=z_j$, $j=1,2,\ldots,N$.
We now examine the possibility of assigning a different set
of collocation nodes. We argue as done in
\cite{miolibro2}, \cite{fatfun} and \cite{funsup}.
\par\smallskip

We start by introducing a function $\chi_N^{\alpha ,\beta }$ as follows:
\begin{equation}\label{fchi}
\chi_N^{\alpha ,\beta } (x)= (1+x)^{\beta }(1-x)P^{\alpha ,\beta }_{N-1}(x) , \quad 
\quad x \in \left[ -1,1 \right],
\end{equation}
and we consider as nodes $x_j$, $j=0,1,\ldots,N$, the zeros  of (\ref{fchi}), that
automatically include the endpoints $\pm 1$. Obviously, $\chi_N^{\alpha ,\beta }$ vanishes at all
points of the grid, thus, the discrete derivative $D_N^\sigma$ applied
to $\chi_N^{\alpha ,\beta }$ is identically zero (viewed from the discrete space,
$\chi_N^{\alpha ,\beta }$ is practically the zero function).
We now apply the exact fractional operator $D^\sigma$ to  $\chi_N^{\alpha ,\beta }$.
This turns out to be an oscillating function.
Successively, we look for collocation nodes such that $[D^\sigma \chi_N^{\alpha ,\beta }](z_i)=0$,
$i=1,2,\ldots,N$. By this choice, we must also get:
\begin{equation}\label{supc}
[(D^\sigma -D_N^\sigma) \chi_N^{\alpha ,\beta }](z_i)=0 ~~~~~~\quad i=1,2,\ldots,N.
\end{equation}
Such an equation tells us that the operator $D^\sigma -D_N^\sigma$ not
only vanishes on the approximation space (by construction, considering
that such a space is the one spanned by the Lagrange type basis (\ref{acca})), but also 
that the extra element $\chi_N^{\alpha ,\beta }$ belongs to the kernel of $D^\sigma -D_N^\sigma$.
This means that we are able to enlarge the dimension of the approximation space
by one unity. As we will check, such an improved consistency property
(called {\sl superconsistency}, according to  \cite{funsup}) is the key
to obtain optimal numerical results, especially when the degree $N$
is not large.
\par\smallskip

In this section, we study both  the ultraspherical case, i.e.:
$\alpha =\beta $   (Subsection \ref{sec41})
 and the case $\alpha =-\beta$  (Subsection \ref{sec42})
although what we are going to say also holds
in a general context. 

%%%%%%%%%%%%%%%%%%%%%%%%%%%%%%%%%%%%%%%%%%%
\subsection{The case $\boldsymbol\alpha =\boldsymbol\beta $}\label{sec41}
%%%%%%%%%%%%%%%%%%%%%%%%%%%%%%%%%%%%%%%%%%%
%
Let us consider the function:
\begin{equation}\label{fchialphabeta}
\chi_N^{\alpha ,\alpha } (x)= (1+x)^{\alpha }(1-x)P^{\alpha ,\alpha }_{N-1}(x) , 
\quad \quad x \in \left[ -1,1 \right],
\end{equation}
and let us assume that the nodes $x_j$, $j=0,1,\ldots,N$, are the zeros of (\ref{fchialphabeta}).

Note that $T^\prime_N$ is proportional to $P^{\alpha ,\beta }_{N-1}$
with $\alpha =\beta =1/2$ (see (\ref{derip}) and (\ref{Cheb})). Thus, the Chebyshev case is included in
our analysis. In the same way, $P^\prime_N$ is proportional to $P^{\alpha ,\beta }_{N-1}$
with $\alpha =\beta =1$  (see (\ref{derip})), so that the Legendre case is also included.
\par\smallskip

The next step is to develop the polynomial $(1-x)P^{\alpha ,\alpha }_{N-1}$ in (\ref{fchialphabeta})
in terms of Jacobi polynomials of the same family. We recall the
recurrence relation  (\ref{recursive2}) when $\alpha=\beta$, so we have:
$$
P^{\alpha ,\alpha }_{N}(x)=a_N xP^{\alpha ,\alpha }_{N-1}(x) 
-c_NP^{\alpha ,\alpha }_{N-2}(x),
$$
\begin{equation}\label{jacr}
{\rm with}~~~~~~~a_N=\frac{(N+\alpha )(2N+2\alpha -1)}{N(N+2\alpha )},
~~~~~c_N=\frac{(N+\alpha )(N+\alpha -1)}{N(N+2\alpha )}.
\end{equation}
Therefore:
\begin{equation}\label{jacre}
\chi_N^{\alpha ,\alpha }(x)= (1+x)^\alpha\left( P^{\alpha ,\alpha }_{N-1}(x)
-\frac{P^{\alpha ,\alpha }_N(x)}{a_N}
-\frac{c_N P^{\alpha ,\alpha }_{N-2}(x)}{a_N}\right).
\end{equation}
To avoid cumbersome calculations, we continue the discussion 
with the case of the Chebyshev Gauss-Lobatto nodes, corresponding
to $\alpha =\beta =1/2$ (generalizations are however straightforward). The above formula 
(\ref{jacre}) takes the form:
$$
\chi_N^{1/2 ,1/2 }(x)= \sqrt{1+x}\left( P^{1/2 ,1/2 }_{N-1}(x)
-\frac{N+1}{2N+1}P^{1/2 ,1/2 }_N(x)
\right.
$$
\begin{equation}\label{jacre2}
\left.
-\frac{2N-1}{4N} P^{1/2 ,1/2 }_{N-2}(x)\right).
\end{equation}

Let $\mu=1-\sigma$. We apply the derivative operator $D^\sigma$
to $\chi_N^{1/2 ,1/2 }$. By plugging the combination (\ref{jacre2}) into (\ref{auto}) 
one obtains:
$$
\Psi_N^{1/2 ,1/2 }(x)=(D^\sigma \chi_N^{1/2 ,1/2 } )(x)=\frac{1}{\Gamma (\mu )}\frac{d}{dx}\int_{-1}^x
\frac{\chi_N^{1/2 ,1/2 }(s) ds}{(x-s)^{1-\mu }}
$$
$$
=\frac{d}{dx}\left[ (1+x)^{1/2+\mu }\left( \frac{\Gamma (N+\frac{1}{2})}
{\Gamma (N+\frac{1}{2} +\mu )}P^{1/2-\mu ,1/2+\mu }_{N-1}(x)\right.\right.
$$
$$
-\frac{N+1}{2N+1}\frac{\Gamma (N+\frac{3}{2})}
{\Gamma (N+\frac{3}{2} +\mu )}P^{1/2-\mu ,1/2+\mu }_N(x)
$$
\begin{equation}\label{fung}
\left.\left. 
 -\frac{2N-1}{4N}\frac{\Gamma (N-\frac{1}{2})}
{\Gamma (N-\frac{1}{2} +\mu )}P^{1/2-\mu ,1/2+\mu }_{N-2}(x)
\right)\right] , 
\end{equation}
where we defined $\Psi_N^{1/2 ,1/2 }$ to be the resulting function.
\par\smallskip

Also easy is the case of the Legendre Gauss-Lobatto nodes, corresponding
to $\alpha =\beta =1$. We have:
\begin{equation}\label{jacre3}
\chi_N^{1,1}(x)= (1+x)\left[ P^{1,1}_{N-1}(x)
-\frac{N(N+2)P^{1,1 }_N(x)}{(N+1)(2N+1)}
-\frac{N P^{1,1 }_{N-2}(x)}{2N+1} \right].
\end{equation}
Successively, one gets:
$$
\Psi_N^{1,1}(x)=(D^\sigma \chi_N^{1,1} )(x)=\frac{1}{\Gamma (\mu )}\frac{d}{dx}\int_{-1}^x
\frac{\chi_N^{1,1}(s) ds}{(x-s)^{1-\mu }}
$$
$$
=\frac{d}{dx}\left[ (1+x)^{1+\mu}\left( \frac{\Gamma (N+1)}
{\Gamma (N+1 +\mu )}P^{1-\mu ,1+\mu }_{N-1}(x)\right.\right.
$$
$$
-\frac{N(N+2)}{(N+1)(2N+1)}\frac{\Gamma (N+2)}
{\Gamma (N+2 +\mu )}P^{1-\mu ,1+\mu }_N(x)
$$
\begin{equation}\label{fung2}
\left.\left. 
 -\frac{N}{2N+1}\frac{\Gamma (N)}
{\Gamma (N +\mu )}P^{1-\mu ,1+\mu }_{N-2}(x)
\right)\right].
\end{equation}
\par\smallskip

\begin{figure}[h!]
\vspace{.1cm}
\centering
{\includegraphics[width=5.8cm,height=4.8cm]{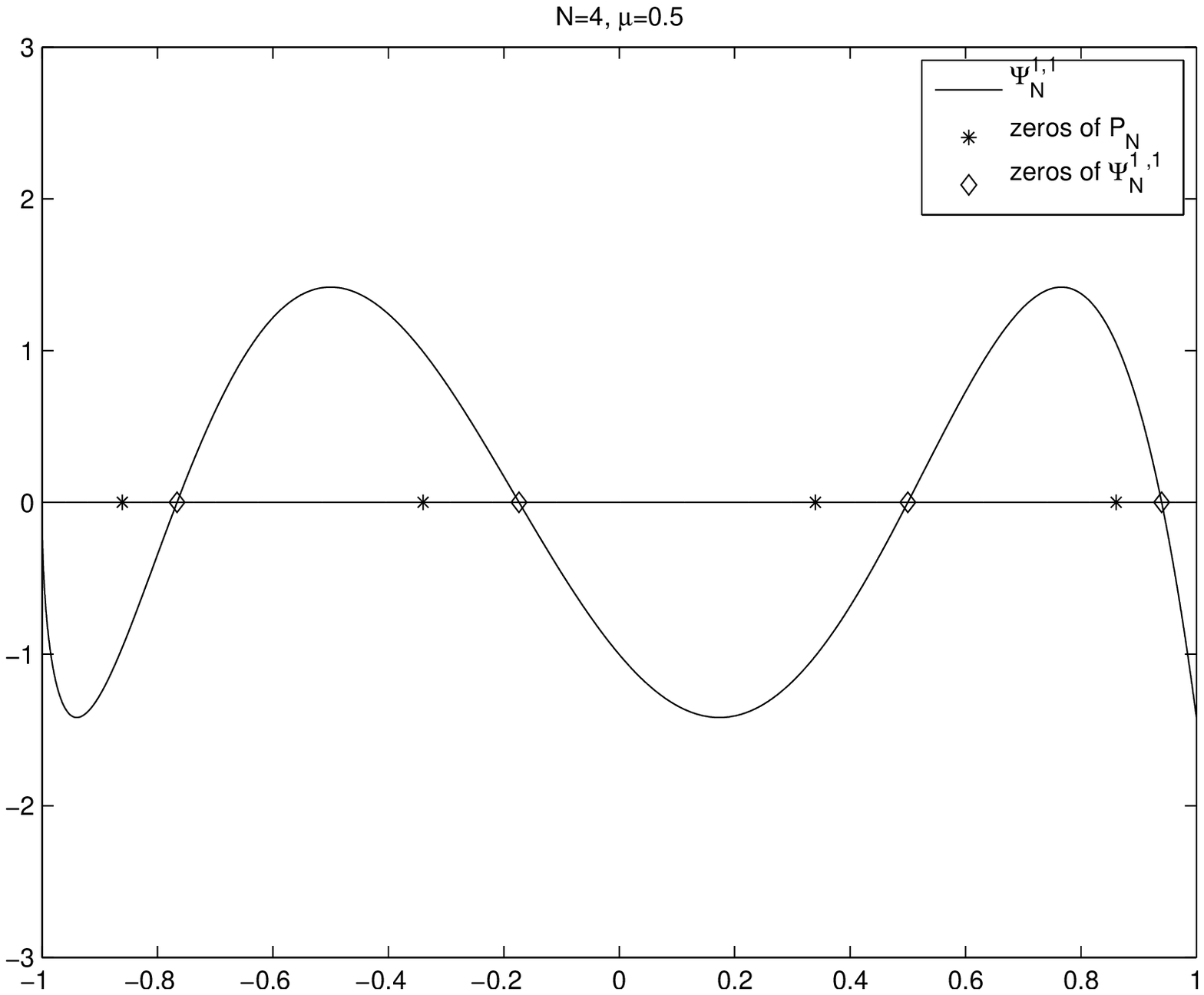}
\hspace{.3cm}\includegraphics[width=5.8cm,height=4.8cm]{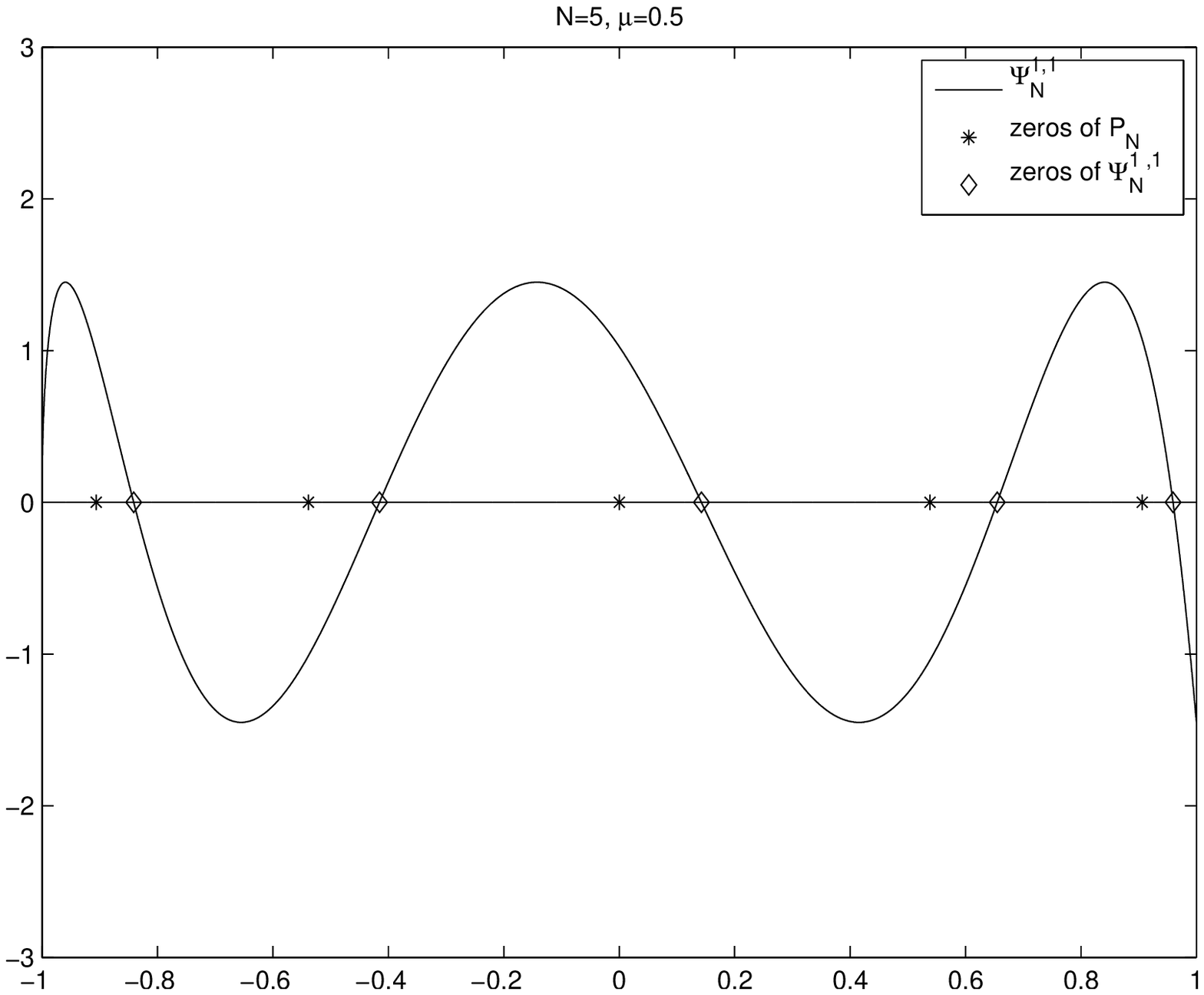}}
\caption{\small \sl  The function  $\Psi^{1 ,1}_N$ for $\mu=0.5$ and  $N=4$ (left), 
$N=5$ (right).}  \label{figPsiN11a}
%\end{figure}
%
%\begin{figure}[h!]
\vspace{.5cm}
\centering
{\includegraphics[width=5.8cm,height=4.8cm]{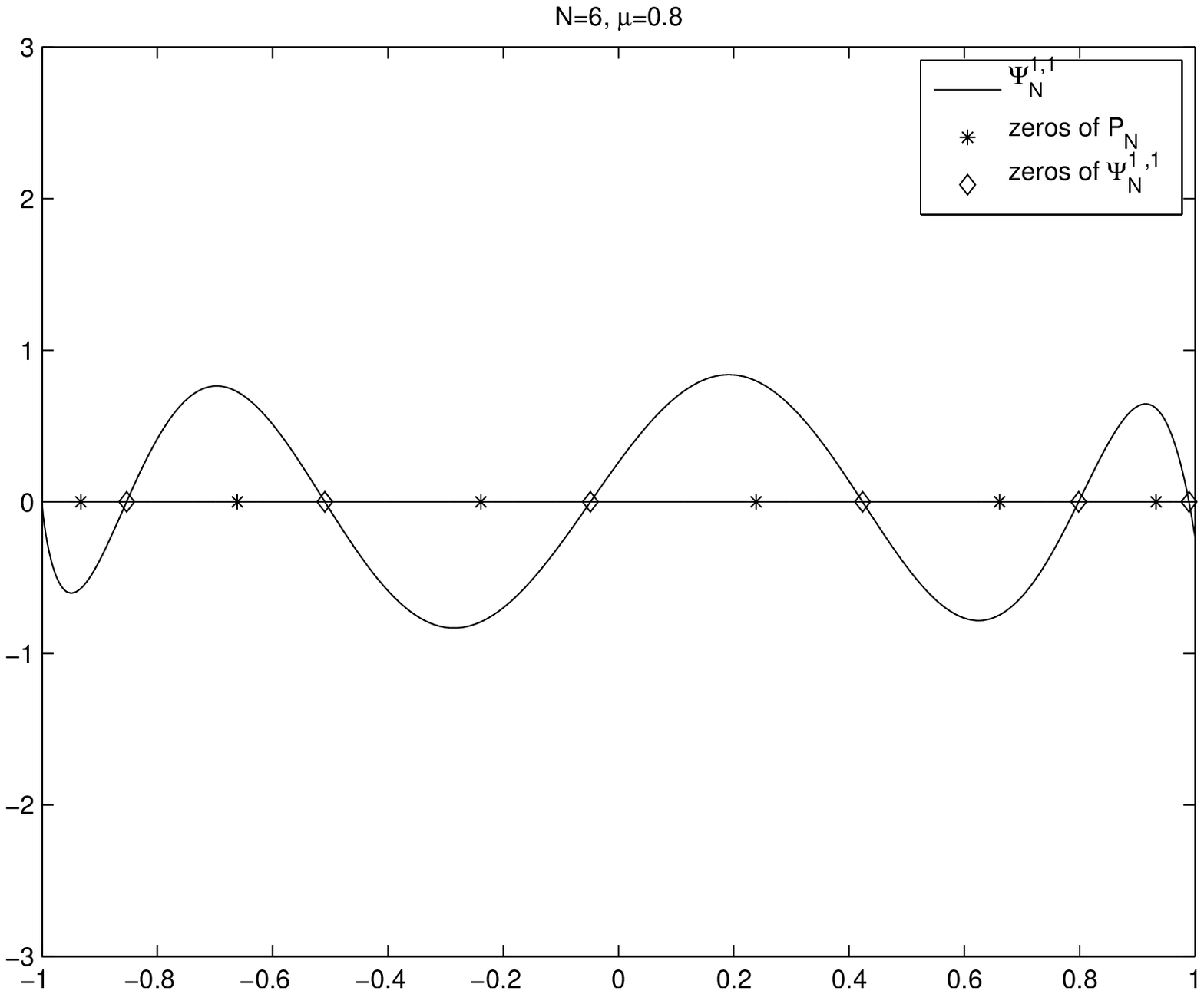}
\hspace{.3cm}\includegraphics[width=5.8cm,height=4.8cm]{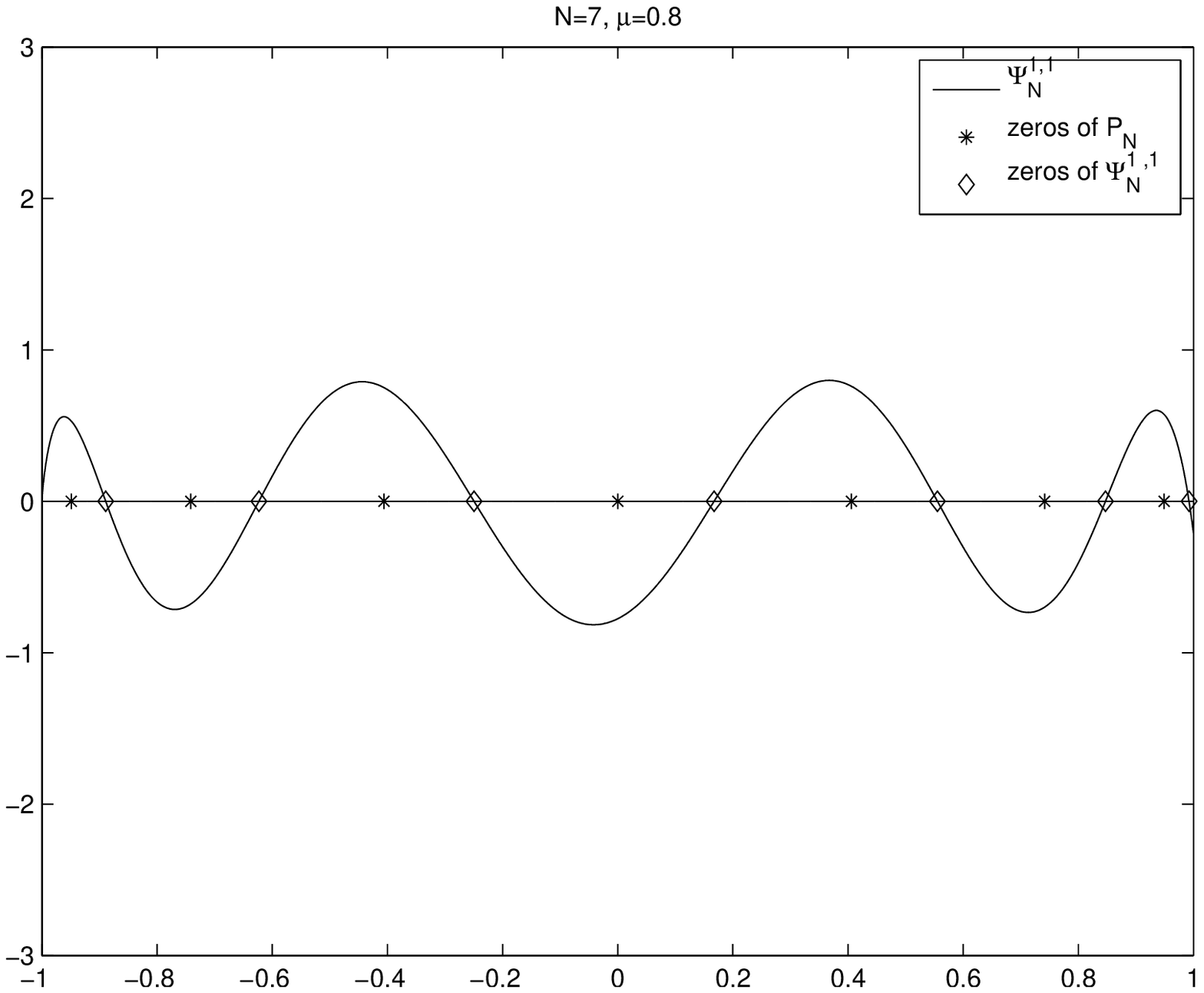}}
\caption{\small \sl  The function  $\Psi^{1 ,1}_N$ for $\mu=0.8$ and  $N=6$ (left), 
$N=7$ (right).}  \label{figPsiN11b}
\end{figure}

We propose to take the collocation nodes $z_i$, $i=1,2,\ldots,N$, to be
the zeros of $\Psi_N^{\alpha ,\beta }$ for $\alpha =\beta =1/2$ or $\alpha =\beta =1$. 
This actually corresponds to find the zeros of a suitable polynomial of degree $N$.
Note that $x=1$ is not a collocation point.
\par\smallskip

For some values of $N$ and $0<\mu <1$, we show in Figures \ref{figPsiN11a}, \ref{figPsiN11b} the plots of the
function  $\Psi_N^{1,1}$. The situation is qualitatively the same for
other values of the parameters. It turns out that we are rather lucky:
there are actually $N$ zeros of $\Psi_N^{1,1}$ in the interval $[-1,1]$ (symbol $\diamond$ in Figures \ref{figPsiN11a}, \ref{figPsiN11b}). Therefore, the
whole construction seems to have solid foundations, though we do not have
a strict proof of this fact. We can point out another singular property. In the pictures we
also plotted the zeros of the Legendre polynomial $P_N$ (symbol $*$ in Figures \ref{figPsiN11a}, \ref{figPsiN11b}). They are interlaced
with our new zeros.  This is quite important for numerical reasons,
if for instance one wants to seek the collocation nodes through the bisection
method.

\par\smallskip

%%%%%%%%%%%%%%%%%%%%%%%%%%%%%%%%%%%%%%%%%%%
\subsection{The case $\boldsymbol\alpha =-\boldsymbol\beta $}\label{sec42}
%%%%%%%%%%%%%%%%%%%%%%%%%%%%%%%%%%%%%%%%%%%

Let us consider the case $\alpha =\mu$, $\beta=-\mu$, and   assume  that the nodes $x_j$, $j=0,1,\ldots,N$, are the zeros of 
$\frac{d}{dx}P^{\mu ,-\mu}_N$ plus the endpoints $\pm 1$. 
As in the previous section, we would like to determine 
the collocation points in order to enlarge the discretization
space by including one more function in the kernel of the
operator $D^\sigma -D^\sigma_N$.
\par\smallskip

Similarly to what has been done before, we  consider the zeros $x_j$, $j=0,1,\ldots,N$, of  the following function:
\begin{equation}\label{fchim}
\chi_N^{1+\mu,1-\mu} (x)= (1+x)^{1-\mu }(1-x)\frac{d}{dx}P^{\mu ,-\mu}_N(x).
\end{equation}
The notation of formula (\ref{fchim}) is justified from the fact
 that $\frac{d}{dx}P^{\mu ,-\mu}_N$
is proportional to $P^{1+\mu  ,1-\mu}_{N-1}$ (see (\ref{derip})), so that:
\begin{equation}\label{fchim2}
\chi_N^{1+\mu,1-\mu} (x)= \frac{N+1}{2} \, (1+x)^{1-\mu }(1-x) \, P^{1+\mu  ,1-\mu}_{N-1}.
\end{equation}

\noindent
From (\ref{recursive2}) we have:
\begin{equation}
P^{1+\mu  ,1-\mu}_{N} (x)=\left(  a_{N} x +b_{N}  \right) P^{1+\mu  ,1-\mu}_{N-1} (x)+c_{N} P^{1+\mu  ,1-\mu}_{N-2} (x), \qquad \forall  n\ge 2 ,
\end{equation}
\smallskip
where:
\begin{eqnarray}
&& a_{N} =\frac{(N+1)(2N+1)}{N(N+2)}    ,    \label{reca1}\\[3mm]
&& b_{N} = \frac{\mu(2N+1)}{N^{2}(N+2)}     ,    \label{recb2}\\[3mm]
&&c_{N} =  - \frac{(N^{2}-\mu^{2})(N+1)}{N^{2}(N+2)}    , \qquad \forall  N\ge 2.\label{recc1}
\end{eqnarray}
As a consequence  the following expression for $\chi_N^{1+\mu,1-\mu}$ holds:
$$
\chi_N^{1+\mu,1-\mu} (x) =(1+x)^{1-\mu }~\frac{N+1}{2}\left[\left(  1+\frac{b_{N}}{a_{N}}\right) P^{1+\mu  ,1-\mu}_{N-1}(x)
\right.
$$
\begin{equation}\label{exppn}
~~~~~~~~~~~~\left. -\frac{1}{a_{N}}P^{1+\mu  ,1-\mu}_N (x)+\frac{c_{N}}{a_{N}}
P^{1+\mu  ,1-\mu}_{N-2}(x)\right].
\end{equation}
On the other hand, by virtue of (\ref{derip}), (\ref{auto}) and (\ref{jac}),
one discovers that:
$$
\int_{-1}^x \frac{(1+s)^{1-\mu } P^{1+\mu  ,1-\mu}_n (s)}
{(x-s)^{1-\mu }} ds=\frac{\Gamma (\mu )\Gamma (n+2-\mu )}{(n+1)!}
(1+x)P^{1,1}_n (x)
$$
\begin{equation}\label{polex}
~~~~~~~~~=\frac{2 \Gamma (\mu )\Gamma (n+2-\mu )}{(n+2)!}
(1+x)P^\prime_{n+1} (x) . 
\end{equation}
Therefore, by applying $D^\sigma$ (with $\mu=1-\sigma$) to the expression of $\chi_N^{1+\mu,1-\mu}$ 
given in (\ref{exppn}),  one finally gets:
$$
\Psi^{1+\mu ,1-\mu}_N(x)=(D^\sigma \chi_N^{1+\mu,1-\mu} )(x)
$$
$$
=(N+1)\frac{d}{dx}\left\{(1+x)\left[\left(  1+\frac{b_{N}}{a_{N}}\right)  \frac{\Gamma (N+1-\mu)}{(N+1)!}
P^\prime_N(x)\right.\right.
$$
$$
\left. \left. -\frac{1}{a_{N}}\frac{\Gamma (N+2-\mu)}{(N+2)!}P^\prime_{N+1}(x)
 +\frac{c_{N}}{a_{N}}\frac{\Gamma (N-\mu )}{ N!}P^\prime_{N-1}(x)
\right]\right\} 
$$
$$
=\frac{(N+1)\Gamma (N-\mu )}{N!}\frac{d}{dx}\left\{(1+x)\frac{d}{dx}
\left[\left(  1+\frac{b_{N}}{a_{N}}\right)  \frac{N-\mu}{N+1}
P_N(x)\right.\right.
$$
\begin{equation}\label{funp}
\left. \left. -\frac{1}{a_{N}}\frac{ (N+1-\mu )(N-\mu )}{(N+2)(N+1)}P_{N+1}(x)
 +\frac{c_{N}}{a_{N}}P_{N-1}(x)
\right]\right\} . 
\end{equation}
\par\smallskip

The right-hand side in (\ref{funp}) is a polynomial of degree $N$.
Thus, we suggest to choose the zeros of $\Psi^{1+\mu ,1-\mu}_N$
as collocation nodes.
\par\smallskip

We can simplify (\ref{funp}) by introducing some approximation.
For $N$ large, one has: $a_N\approx 2$, $c_N\approx -1$, while
$b_N$ tends to zero. Thus, for  $N$ large we are allowed to write:
$$
\frac{d}{dx}\left\{(1+x)\frac{d}{dx}
\left[\left(  1+\frac{b_{N}}{a_{N}}\right)  \frac{N-\mu}{N+1}
P_N(x)\right.\right.
$$
$$
\left. \left. -\frac{1}{a_{N}}\frac{ (N+1-\mu )(N-\mu )}{(N+2)(N+1)}P_{N+1}(x)
 +\frac{c_{N}}{a_{N}}P_{N-1}(x)
\right]\right\}  
$$
\begin{equation}\label{funp2}
\approx \frac{d}{dx}\left\{(1+x)\frac{d}{dx}
\left[ P_N(x)-\frac{1}{2}P_{N+1}(x)-\frac{1}{2}P_{N-1}(x)
\right]\right\} . 
\end{equation}
On the other hand, by the recurrence relation for Legendre polynomials
(see (\ref{recursive2}) for $\alpha =\beta =0$), we deduce that:
\begin{equation}\label{legre}
\frac{1}{2}P_{N+1}(x)+\frac{1}{2}P_{N-1}(x)\approx x P_N(x).
\end{equation}
In this way, the right-hand side of (\ref{funp2}) is approximated by:
$$
\frac{d}{dx}\left\{(1+x)\frac{d}{dx}
\left[ (1-x)P_N
\right]\right\} 
$$
\begin{equation}\label{funp3}
=\Big( (1-x^2)P^\prime_N\Big)^\prime
-\Big( (1+x)P_N \Big)^\prime =-(N^2+N+1)P_N- (1+x)P_N^\prime , 
\end{equation}
where we used the differential equation characterizing Legendre
polynomials (see (\ref{slp}) for $\alpha =\beta =0$). It is worthwhile to observe 
that these last formulas do not depend on $\mu$.
\par\smallskip

Thus, from (\ref{funp})  for $N$ sufficiently large we can write:
\begin{equation}\label{funpapprox}
\Psi^{1+\mu ,1-\mu}_N\approx -\frac{(N+1)\Gamma (N-\mu )}{N!}\left[  
(1+x)P_N^\prime+(N^2+N+1)P_N \right],
\end{equation} 
that could turn out to be useful in understanding the theoretical properties
of the function $\Psi^{1+\mu ,1-\mu}_N$, such as the location of
its zeros.
Note that in the above formula the multiplying constant approaches $(N+1)/N$
for $\mu \rightarrow 0$ and $(N+1)/N(N-1)$ for $\mu \rightarrow 1$.
\par\smallskip

\begin{figure}[h!]
\vspace{.1cm}
\centering
{\includegraphics[width=5.8cm,height=4.8cm]{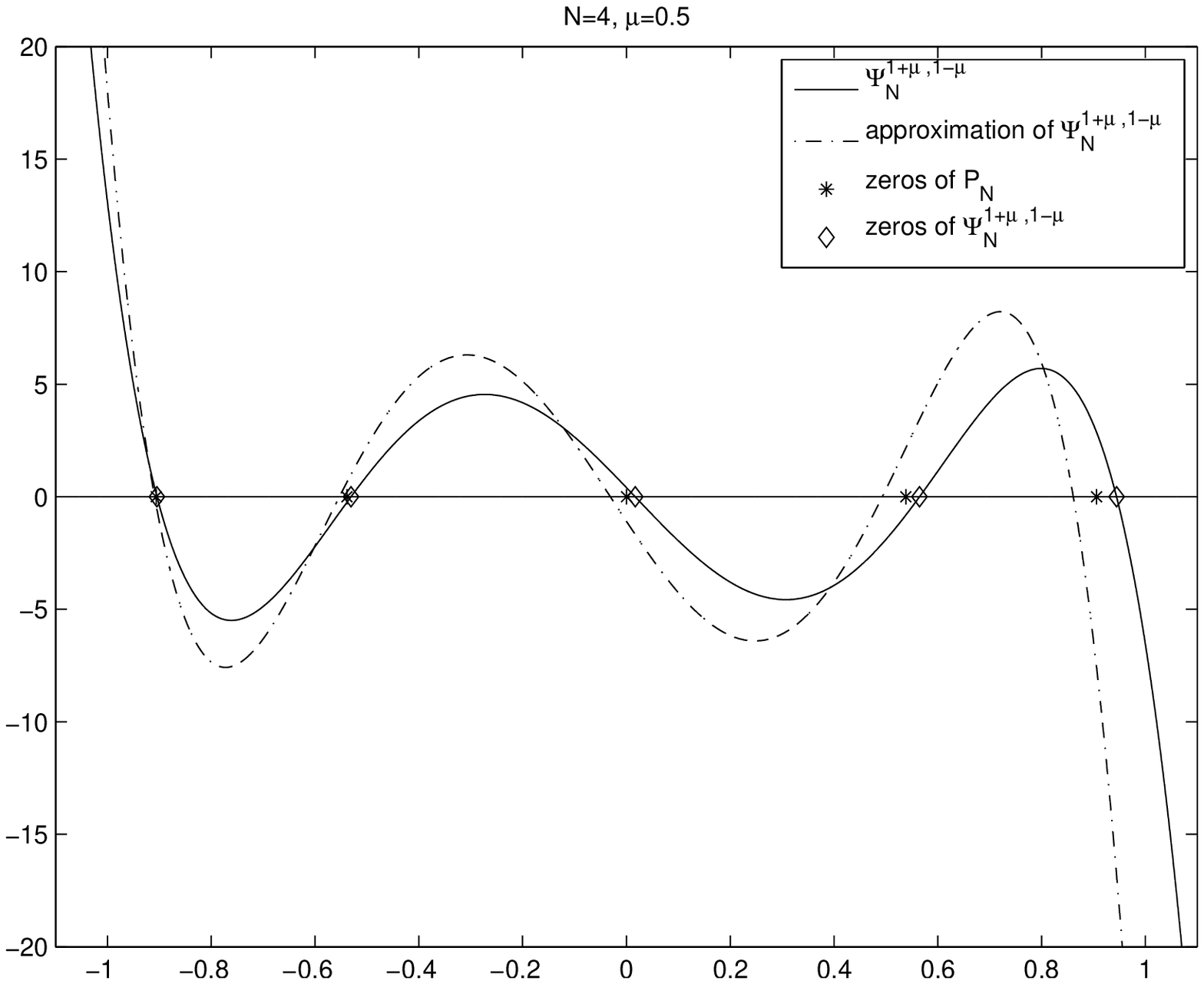}
\hspace{.3cm}\includegraphics[width=5.8cm,height=4.8cm]{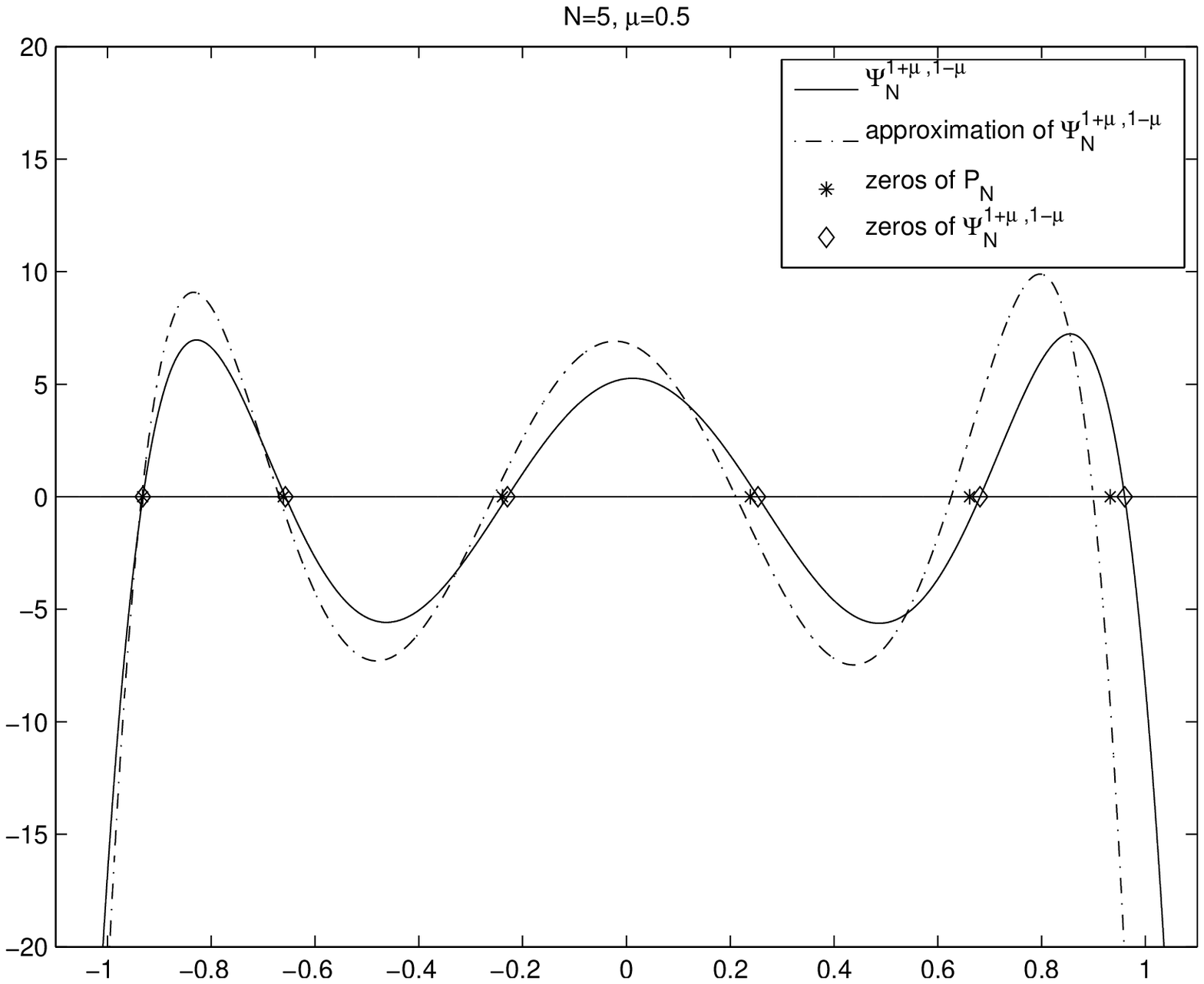}}
\caption{\small \sl  The function  $\Psi^{1+\mu ,1-\mu}_N$, $\mu=0.5$,  and its approximation for $N=4$ (left)  and
$N=5$ (right).}  \label{fig3}
\end{figure}

We expect that the zeros of $\Psi^{1+\mu ,1-\mu}_N$
are not very far from those of the right-hand side in  (\ref{funpapprox}).
As a matter of fact, we compare in Figure \ref{fig3}  the plots of $\Psi^{1+\mu ,1-\mu}_N$
with those of the corresponding approximations, for different values of $N$ and $\mu$.
In both cases, 
there are exactly $N$ zeros, all belonging to the interval $]-1,1[$. A comparison with the 
zeros of $P_N$ is also made (symbol $*$ in Figure \ref{fig3}). 
The two sets of points (zeros of $\Psi^{1+\mu ,1-\mu}_N$ and Legendre zeros) are
alternate. This has proven to be true for all the values of $\mu$ and $N$ that we tested.
Having an idea of the distribution of the new collocation points 
(symbol $\diamond$ in Figure \ref{fig3}) is important in view of developing
methods for their numerical computation. For instance, the implementation of the
bisection methods follows naturally. In order to set up the collocation scheme
that is discussed in the coming sections,
we actually computed the approximate zeros of $\Psi^{1+\mu ,1-\mu}_N$ by the bisection method. 
\par\smallskip

\begin{figure}[h!]
\vspace{.1cm}
\centering
\includegraphics[width=8.cm,height=6.5cm]{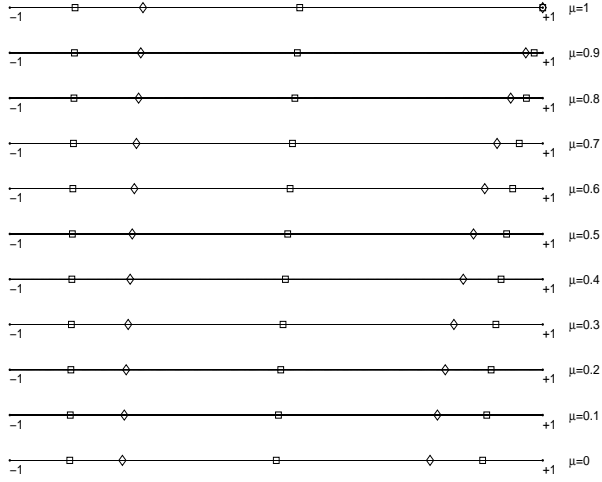}
\caption{\small \sl Zeros of $\Psi^{1+\mu ,1-\mu}_N$ for $N=2$ (symbol $\diamond$) and $N=3$ 
(symbol $\square$) for different values of $\mu$.}  \label{fig5}
\end{figure}

For convenience, in Figure \ref{fig5} we also show the zeros of $\Psi^{1+\mu ,1-\mu}_N$ 
for $N=2$ (symbol $\diamond$) 
and $N=3$ (symbol $\square$) when $\mu$ varies from 0 to 1 with step 0.1.
\par\smallskip

%%%%%%%%%%%%%%%%%%%%%%%%%%%%%%%%%%%%%%%%%%%
\section{Some preliminary numerical results}\label{sec5}
%%%%%%%%%%%%%%%%%%%%%%%%%%%%%%%%%%%%%%%%%%%

Based on what obtained in the previous section, there are
many possible ways to approach the fractional  differential 
problem (\ref{pru}), through  the discrete problem (\ref{prud}), depending on 
the construction of the discrete fractional operator $D^\sigma_{N}$, $0 <\sigma <1$.
In the examples we are discussing, variants do not only rely on the choice of the
initial polynomial basis, but also on the choice of the collocation grid.
\par\smallskip 

In all the cases we are going to examine in this section, given $0 <\sigma <1$, we
take $\mu =1-\sigma$ and  we represent the function
to be differentiated using the grid-points $x_j$, taken to be the $N-1$ zeros
of $\frac{d}{dx}P^{\mu ,-\mu}_N$ plus the endpoints $x_N=\pm 1$. Those correspond to the zeros of  the function 
$\chi_N^{1+\mu,1-\mu} $ defined in  (\ref{fchim}).
One then builds  the $N\times N$ approximation matrices $D^\sigma_{N}$
by suitably choosing the collocation nodes  $z_i$,  $i=1,\ldots,N$.
For simplicity, we just consider the following possibilities:

\begin{description}
\item[ {\sl Choice 1}] - The representation nodes $x_j$, $j=1,2,\ldots,N$,
are the zeros of $\frac{d}{dx}P^{\mu ,-\mu}_N$ plus the endpoint $x_N=1$  and the collocation nodes coincide with the representation nodes, i.e.,  $z_i=x_{i}$, $i=1,2,\ldots,N$;
\item[ {\sl Choice 2}] - The representation   nodes $x_j$, $j=1,2,\ldots,N$,
are the zeros of $\frac{d}{dx}P^{\mu ,-\mu}_N$ plus the endpoint $x_N=1$ and the collocation nodes $z_i$, $i=1,2,\ldots,N$,
are the zeros of the derivative  of the Chebyshev 
polynomial $T_N$  with the addition of the point $x_N=1$, i.e., the points  defined in (\ref{nodiCheb});
\item[ {\sl Choice 3 (superconsistency)}] - The representation   nodes  $x_j$, $j=1,2,.$ $..,N$,
are the zeros of $\frac{d}{dx}P^{\mu ,-\mu}_N$ plus the endpoint $x_N=1$ and the collocation nodes $z_i$, $i=1,2,\ldots,N$,
are the $N$ zeros of  $\Psi^{1+\mu ,1-\mu}_N=D^\sigma \chi_N^{1+\mu,1-\mu}$, as
defined in (\ref{funp}).
 \end{description}
 
The fact that the representation nodes are all the same in these examples will allow us
to make appropriate and consistent comparisons.
\par\smallskip

\begin{table}{h!}
\centerline{
{\begin{tabular}{@{}ccclccclccclccclccc}\hline
 &&$N$ &&  Error 1&&      Error 2 &&  Error 3 &&           \\
\hline
 &&    4   &&       0.4057     &&     0.6852   &&      0.0824    &&     \\
 &&    5   &&       0.2053     &&     0.3807   &&      0.0363   &&      \\
 &&    6   &&       0.1348     &&     0.4069   &&      0.0084   &&      \\
 &&    7   &&       0.0764     &&     0.1365   &&      0.0039   &&      \\
 &&    8   &&       0.0140     &&     0.0316  &&       0.0015   &&      \\
 &&    9   &&       0.0143     &&     0.0417  &&       3.7184e-04 &&      \\
&&    10   &&    0.00653     &&     0.0128  &&      1.1357e-04 \\
&&    11   &&   5.4582e-04 &&     0.0012 &&    4.7035e-05 \\
&&    12   &&   7.0567e-04 &&    0.0021 &&    1.1296e-05 \\
&&    13   &&   3.1975e-04 &&    6.8793e-04 && 1.9664e-06 \\
&&    14   &&   3.7780e-05 &&    2.4258e-05 && 9.9621e-07 \\
&&    15   &&   1.9172e-05 &&    6.1681e-05&& 2.5710e-07\\
 \hline
  \end{tabular}}
  }\vspace{-0.5cm}
  \bigskip\caption{\small \sl Errors in 
the discrete maximum norm between the exact solution  $u$ and the approximated solution $u_{N}$
of problem (\ref{pru}) with $\sigma=0.5$, $g(x)=\sin 2(x+1)^2$. They are  obtained with
the same representation nodes, but with different collocation nodes, as a consequence of
{\sl Choice 1} (Error 1), {\sl Choice 2} (Error 2), and {\sl Choice 3} (Error 3).}
\label{table2}
\end{table}

We solve numerically the fractional differential  problem (\ref{pru}) with $\sigma=0.5$ and
right-hand side $g(x)=\sin 2(x+1)^2$. We examine the three possibilities itemized above.
Note that the exact solution of $D^\sigma u=g$ is not known; therefore,
we compute an approximation of $u$ with $N$ relatively large, to be used 
in place of $u$ in our comparisons. It does not matter what set of collocation
nodes is utilized in this operation, since, due to spectral convergence,
the various approximations are graphically indistinguishable, when $N$
is sufficiently large.  The solution  $u_N$ of (\ref{prud}) is represented through
(\ref{svil}),  with  $H_j$, $j=1,2,\ldots,N$, defined in  (\ref{acca}).

\par\smallskip

\begin{figure}[h!]
\vspace{.1cm}
\centering{\includegraphics[width=5.8cm,height=4.8cm]{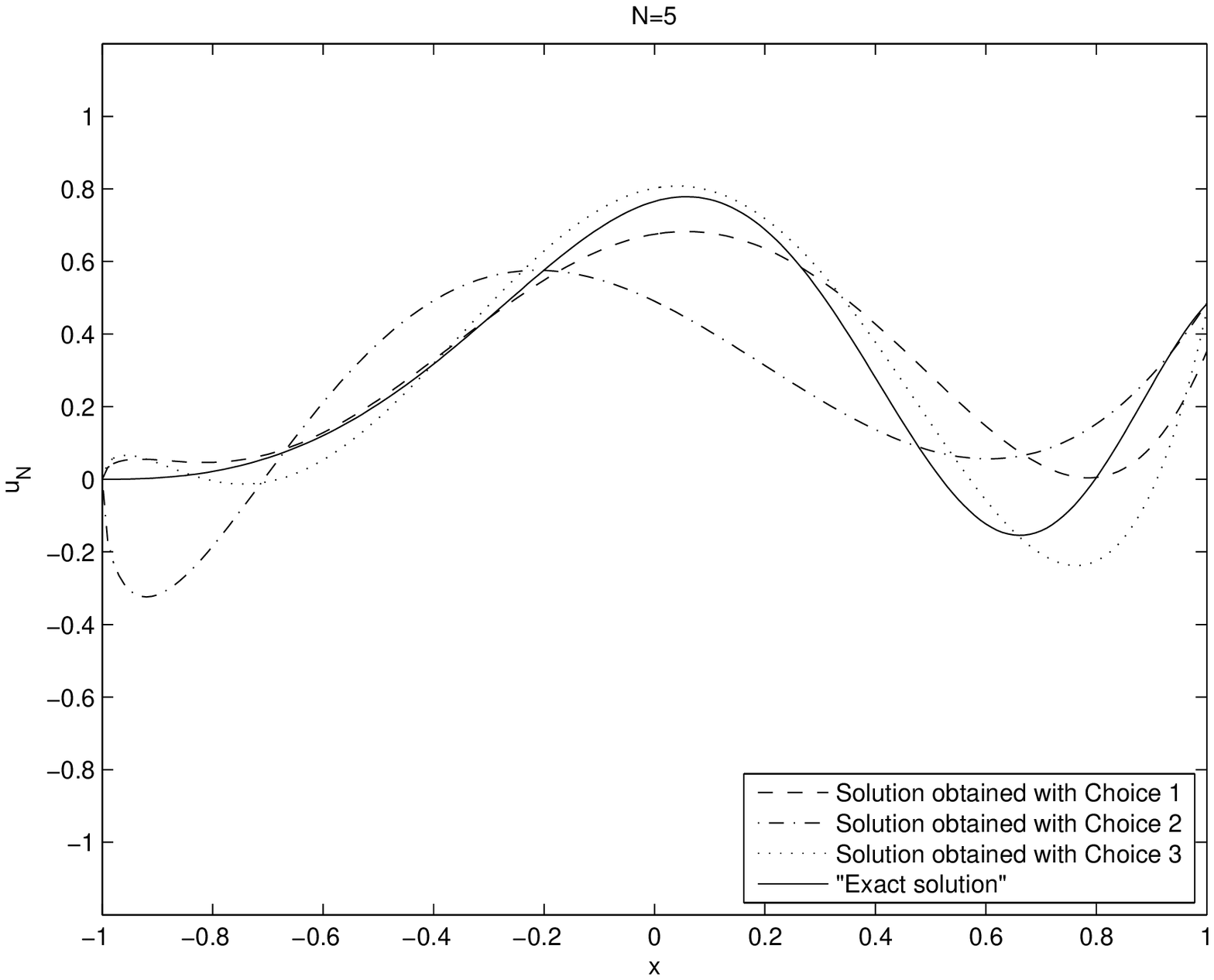}
\hspace{.3cm}\includegraphics[width=5.8cm,height=4.8cm]{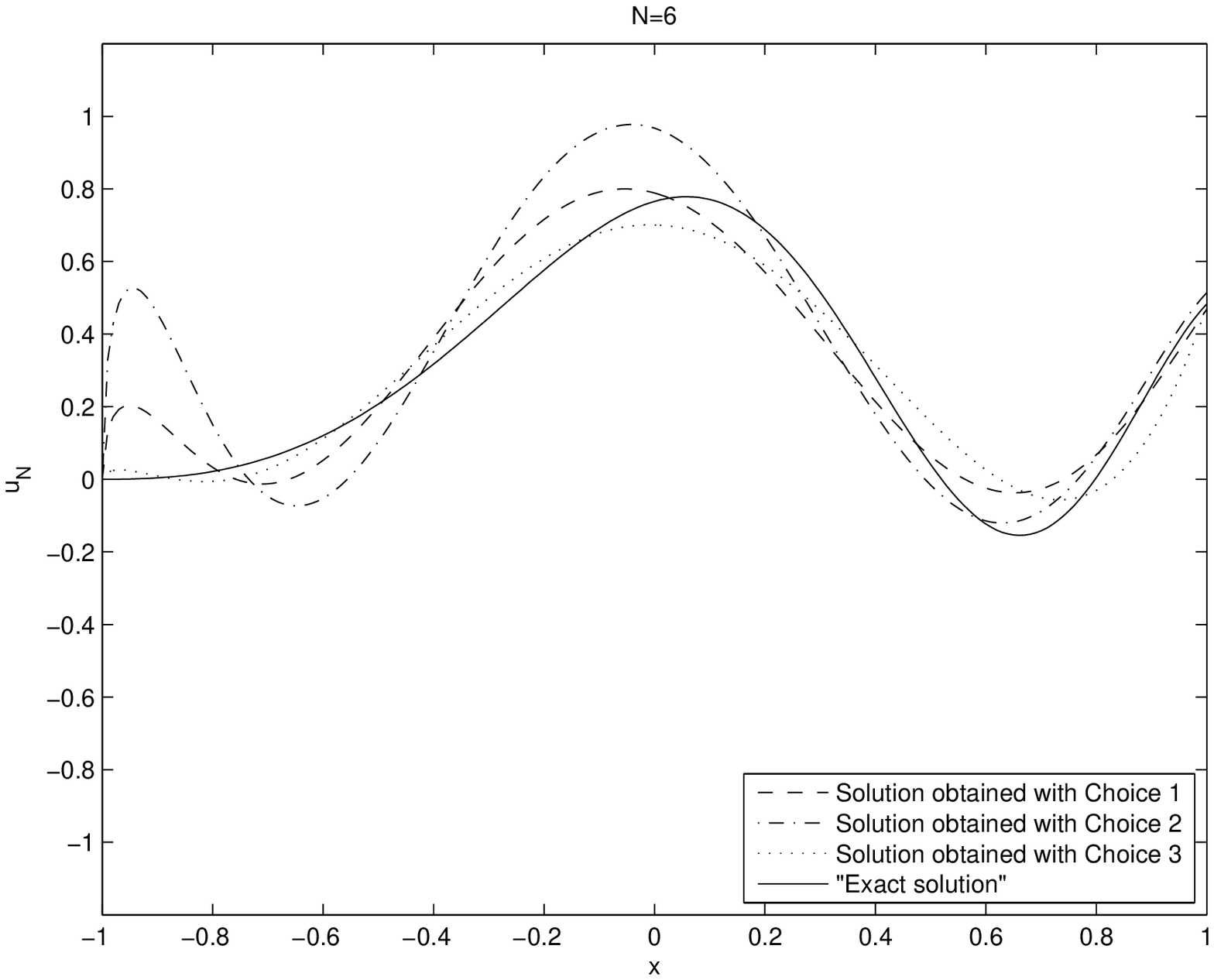}}
\caption{\small \sl 
Approximated solutions for $N=5$ and $N=6$ of the fractional differential  problem (\ref{pru}) 
with $g(x)=\sin 2(x+1)^2$, $\sigma=0.5$, for the three  different choices of  collocation nodes  described  
in Section \ref{sec5}.
}  \label{fig6}
%\end{figure}
%
%\begin{figure}[h!]
\vspace{.1cm}
\centering{\includegraphics[width=5.8cm,height=4.8cm]{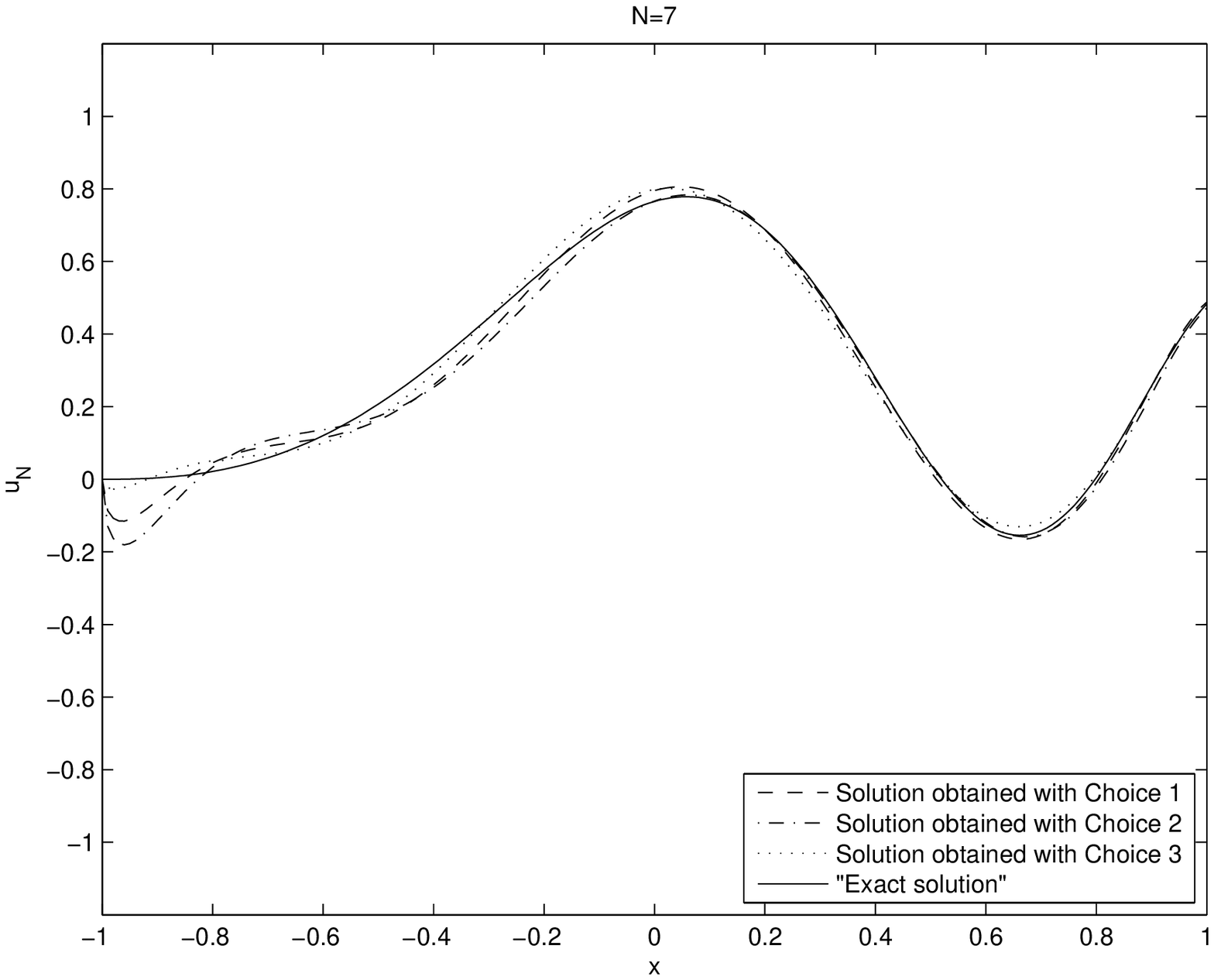}
\hspace{.3cm}\includegraphics[width=5.8cm,height=4.8cm]{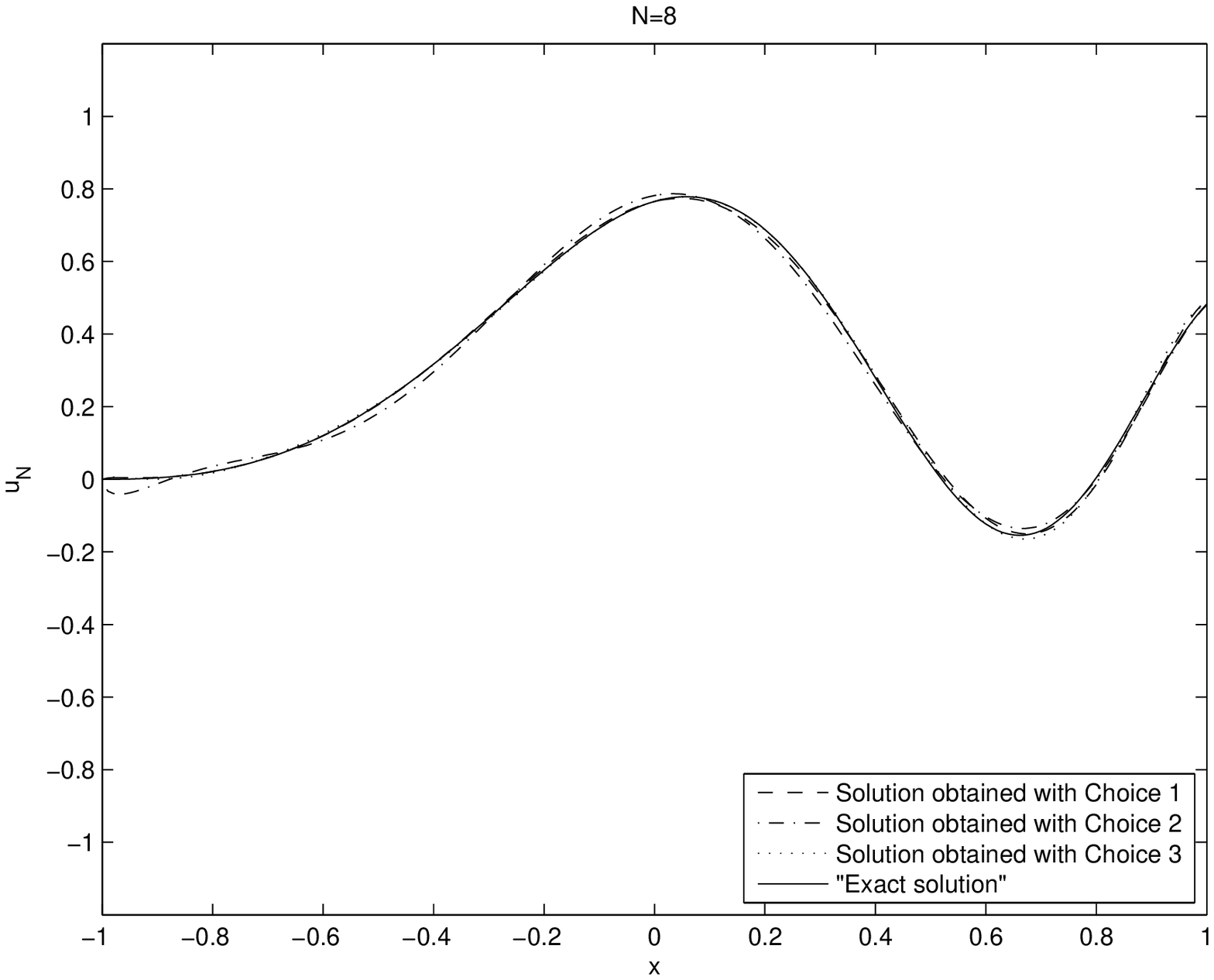}}
\caption{\small \sl 
Approximated solutions for $N=7$ and $N=8$ of the fractional differential  problem (\ref{pru}) 
with $g(x)=\sin 2(x+1)^2$, $\sigma=0.5$, for the three  different choices of  collocation nodes  described  
in Section \ref{sec5}.
}  \label{fig7}
\end{figure}

The results of these tests are given in Figures \ref{fig6} and  \ref{fig7} for $N=5,6,7,8$, 
respectively. Here $N$ is not large,
so that some differences soon emerge, depending on the choice of the collocation sets.   
The ``exact solution'' in Figures \ref{fig6} and \ref{fig7} has been computed with $N=50$ 
using for both representation and collocation nodes the points  defined in (\ref{nodiCheb}).
It is clear from the pictures that the approach here proposed provides reasonable
approximations, even at such low polynomial degrees, while the more classical methods
tend to be less accurate. 
\par\smallskip

By enlarging $N$, all the three types of approximated solutions
converge spectrally. Thus, their plots are almost coincident. By the way, we give
in Table \ref{table2} the errors, in the discrete maximum norm, relative to the cases examined.
The best performances are provided by  {\sl Choice 3} (Error 3), corresponding to the superconsistent
method. The other methods look more erratic; in practice one needs to take $N$ suitably
large before observing a correct decay rate.

\par\bigskip

We end this sections by showing how to recover the entries
of the matrix $A_N$ in (\ref{mat}) in explicit way,
when the representation nodes are the zeros of $\frac{d}{dx}P^{\mu ,-\mu}_N(x_j)$
plus the endpoints $\pm 1$. We start by recalling the Gauss-Lobatto integration formula:
\begin{equation}\label{gau}
\int_{-1}^1q(x)(1-x)^\alpha (1+x)^\beta dx=\sum_{m=0}^N q(x_m)w_m , 
\end{equation}
valid for any $q$ polynomial of degree less or equal to $2N-1$.
The weights $w_m$, $m=1,2,\ldots,N$ are known (see \cite{miolibro}, p. 52). 
In particular, if $\alpha =\mu$
and $\beta =-\mu$, we have for $1\leq m\leq N-1$:
\begin{equation}\label{pesi}
w_m=\frac{2 ~\Gamma (N+\mu ) \Gamma (N-\mu ) }{(N+1)~\left[(N-1)!\right]^2}
~\frac{-1}{P^{\mu ,-\mu}_N(x_m) \frac{d}{dx}P^{\mu ,-\mu}_{N-1}(x_m)}.
\end{equation}
Similar formulas are available for $w_0$ and $w_N$. 

For $1\leq k\leq N-1$, we multiply (\ref{svilhm})  by $P^{\mu ,-\mu}_k(x_m)w_m$;
then we sum up on the index $m$, obtaining: 
$$
\sum_{n=1}^N c(j,n) \sum_{m=0}^N  \left[ P_n^{\mu ,-\mu}(x_m)
-P_n^{\mu ,-\mu}(-1)\right] P_k^{\mu ,-\mu}(x_m)w_m
$$
\begin{equation}\label{svisu}
= \sum_{m=0}^N \delta_{jm}(x_m+1)^{\mu}P_k^{\mu ,-\mu}(x_m)w_m
= (x_j+1)^{\mu}P_k^{\mu ,-\mu}(x_j)w_j . 
\end{equation}
Since $k+n\leq 2N-1$, we can apply Gaussian integration, and
due to the orthogonality of Jacobi polynomials one finally gets:
\begin{equation}\label{svisuo}
c(j,k)\int_{-1}^1 \left[ P_k^{\mu ,-\mu}(x)\right]^2(1-x)^\mu (1+x)^{-\mu}dx
= (x_j+1)^{\mu}P_k^{\mu ,-\mu}(x_j)w_j , 
\end{equation}
from which one easily gets the coefficients for $1\leq j\leq N$
and $1\leq k\leq N-1$. If $k=N$, we can still arrive at the expression
in (\ref{svisu}). Successively, we can use the orthogonality
when the index $n$ is between 1 and $N-1$. In the end, we get:
\begin{equation}\label{svisuon}
c(j,N)\sum_{m=0}^N \left[ P_N^{\mu ,-\mu}(x_m)\right]^2w_m
= (x_j+1)^{\mu}P_N^{\mu ,-\mu}(x_j)w_j ,
\end{equation}
from which one recovers $c(j,N)$, for $1\leq j\leq N$.
\par\smallskip

%
%%%%%%%%%%%%%%%%%%%%%%%%%%%%%%%%%%%%%%%%%%%
\section{Application to a boundary-value problem}\label{sec6}
%%%%%%%%%%%%%%%%%%%%%%%%%%%%%%%%%%%%%%%%%%%
In this last section, we would like to approximate the
solution $u$ of the following differential fractional equation
with homogeneous Dirichlet boundary constraints:
\begin{equation}\label{dife}
-\frac{d^2u}{dx^2} +K D^\sigma u =g,   ~~~~~~~~\quad {\rm with}~u(\pm 1)=0 , 
\end{equation}
where $0<\sigma <1$, and $K$ is a given constant.  Throughout this section we take $\mu =1-\sigma$.
\par\smallskip

In finite dimension,  problem (\ref{difen})  reads as follows:
\begin{equation}\label{difen}
-D^2_Nu_N +KD^\sigma_N u_N =g_N   ~~~~~~~~~~~~~{\rm with}~u_N(\pm 1)=0 , 
\end{equation}
and this equation must hold at some collocation points.
\par\smallskip

The solution  $u_N$ of (\ref{difen}) is represented as in (\ref{svil})  
with  $H_j$, $j=1,2,\ldots,N$, defined in  (\ref{acca}).
Since we want to impose
boundary conditions, the sum in (\ref{svil}) goes from $j=1$
up to $j=N-1$. 
 Once the representation  points $x_j$, $j=1,2,\ldots,N-1$, are chosen, the
approximating matrix is recovered by applying the discrete operator 
$-D^2_N+KD^\sigma_N$ to $H_j$ and evaluating at the collocation nodes $z_i$, $i=1,2,\ldots,N-1$.
\par\smallskip

From the results of the previous sections, several possibilities
may be taken into account both for the representation grid and the collocation grid.
In order to get {\sl superconsistent} type approximations,
the following possibilities are taken into account.

%*********************************
\par\smallskip

\parindent=.7cm {\sl 1)} The representation nodes  $x_j$, $j=1,2,\ldots,N-1$,
are the zeros of $T^\prime_N$ and the collocation nodes $z_i$, $i=1,2,\ldots,N-1$, are chosen  such that: 
\begin{equation}\label{collo1}
-\chi_N^{\prime\prime}(z_i)+K\Psi_{N}(z_i)=0 , 
\end{equation}
where $\chi_N=\chi_N^{1/2,1/2 } $ is defined in (\ref{jacre2}) 
and $\Psi_N=\Psi_N^{1/2 ,1/2 }$ is given in (\ref{fung}). The above
equation actually admits $N-1$ roots in the interval $]-1,1[$. 
This leads us to a squared matrix
of dimension $(N-1)\times (N-1)$. 
\par\medskip

\parindent=.7cm {\sl 2)} The representation nodes $x_j$, $j=1,2,\ldots,N-1$, are the zeros of $P^\prime_N$ and 
the collocation nodes $z_i$, $i=1,2,\ldots,N-1$,  are the solution of  equation (\ref{collo1}) with
 $\chi_N=\chi_N^{1,1 } $  defined in (\ref{jacre3})  
 and $\Psi_N=\Psi_N^{1,1 }$  obtained from (\ref{fung2}).
\par\medskip

\parindent=.7cm {\sl 3)} The
representation nodes $x_j$, $j=1,2,\ldots,N-1$, are the zeros of  $\frac{d}{dx}P^{\mu ,-\mu}_N$
and  the collocation nodes $z_i$, $i=1,2,\ldots,N-1$,  are the solution of  equation (\ref{collo1}) with
  $\chi_N=\chi_N^{1+\mu,1-\mu}$   defined in (\ref{fchim})  
and $\Psi_N=\Psi_N^{1+\mu,1-\mu}$  computed in (\ref{funp}).
\par\medskip

We now concentrate our attention on point 3) by carrying out some tests, aimed
to compare (as we did in Section 5) the technique here proposed with 
the more standard ones. We solve numerically the fractional differential problem (\ref{dife})
in the following circumstances:

\begin{description}
\item[ {\sl Choice 4}] - The representation  points $x_j$, $j=1,2,\ldots,N-1$, are the zeros of  $\frac{d}{dx}P^{\mu ,-\mu}_N$ 
 and the collocation nodes coincide with the representation nodes, i.e.,  $z_i=x_{i}$, $i=1,2,\ldots,N-1$;
\item[ {\sl Choice 5}] - The representation  points $x_j$, $j=1,2,\ldots,N-1$, are the zeros of  $\frac{d}{dx}P^{\mu ,-\mu}_N$ 
 and the collocation nodes $z_i$, $i=1,2,\ldots,N-1$, are the zeros of the derivative  of the Chebyshev 
polynomial $T_N$  plus the point $x_N=1$, i.e. the points  defined in (\ref{nodiCheb});
\item[ {\sl Choice 6 }] ({\sl Superconsistency}, see option 3) above) - 
The representation  points $x_j$, $j=1,2,\ldots,N-1$, are the zeros of  $\frac{d}{dx}P^{\mu ,-\mu}_N$ 
 and  the collocation nodes $z_i$, $i=1,2,\ldots,N-1$,
are the $N-1$ zeros of (\ref{collo1})  where $\chi_N=\chi_N^{1+\mu,1-\mu}$  is defined in (\ref{fchim})  
and $\Psi_N=\Psi_N^{1+\mu,1-\mu}$ is computed in (\ref{funp}).

 \end{description}

Note that in order to approach the boundary-value problem  (\ref{dife})  
we also need to evaluate the second derivative of $\chi_N=\chi_N^{1+\mu,1-\mu}$  defined in (\ref{fchim}).  
To this scope, let us note that  from  (\ref{slp}) with  
$\alpha =\mu$ and $\beta =-\mu$, one has:
\begin{equation}\label{difej}
\frac{d}{dx}\left[ (1-x^2)\frac{d}{dx}P^{\mu ,-\mu}_N\right]
=2\mu\frac{d}{dx}P^{\mu ,-\mu}_N-N(N+1)P^{\mu ,-\mu}_N, 
\end{equation}
and consequently:
$$
\frac{d^2}{dx^2} \left[ (1+x)^{-\mu}(1-x^2) \frac{d}{dx}P^{\mu ,-\mu}_N\right]
=\mu (\mu +1) (1+x)^{-\mu -2}(1-x^2) \frac{d}{dx}P^{\mu ,-\mu}_N
$$
$$
-2\mu (1+x)^{-\mu -1}\frac{d}{dx} \left[ (1-x^2) \frac{d}{dx}P^{\mu ,-\mu}_N\right]
+(1+x)^{-\mu}\frac{d^2}{dx^2} \left[ (1-x^2) \frac{d}{dx}P^{\mu ,-\mu}_N\right]
$$
$$
=(1+x)^{-\mu -2}\left[ \mu (\mu +1)  (1-x^2)-4\mu^2 (1+x)
-N(N+1)(1+x)^2\right] \frac{d}{dx}P^{\mu ,-\mu}_N
$$
\begin{equation}\label{chisec}
+2\mu  N(N+1)(1+x)^{-\mu -1}P^{\mu ,-\mu}_N +2\mu(1+x)^{-\mu }\frac{d^2}{dx^2}P^{\mu ,-\mu}_N.
\end{equation}
\par\smallskip

\begin{figure}[h!]
\vspace{.1cm}
\centering
\centering{\includegraphics[width=5.8cm,height=4.8cm]{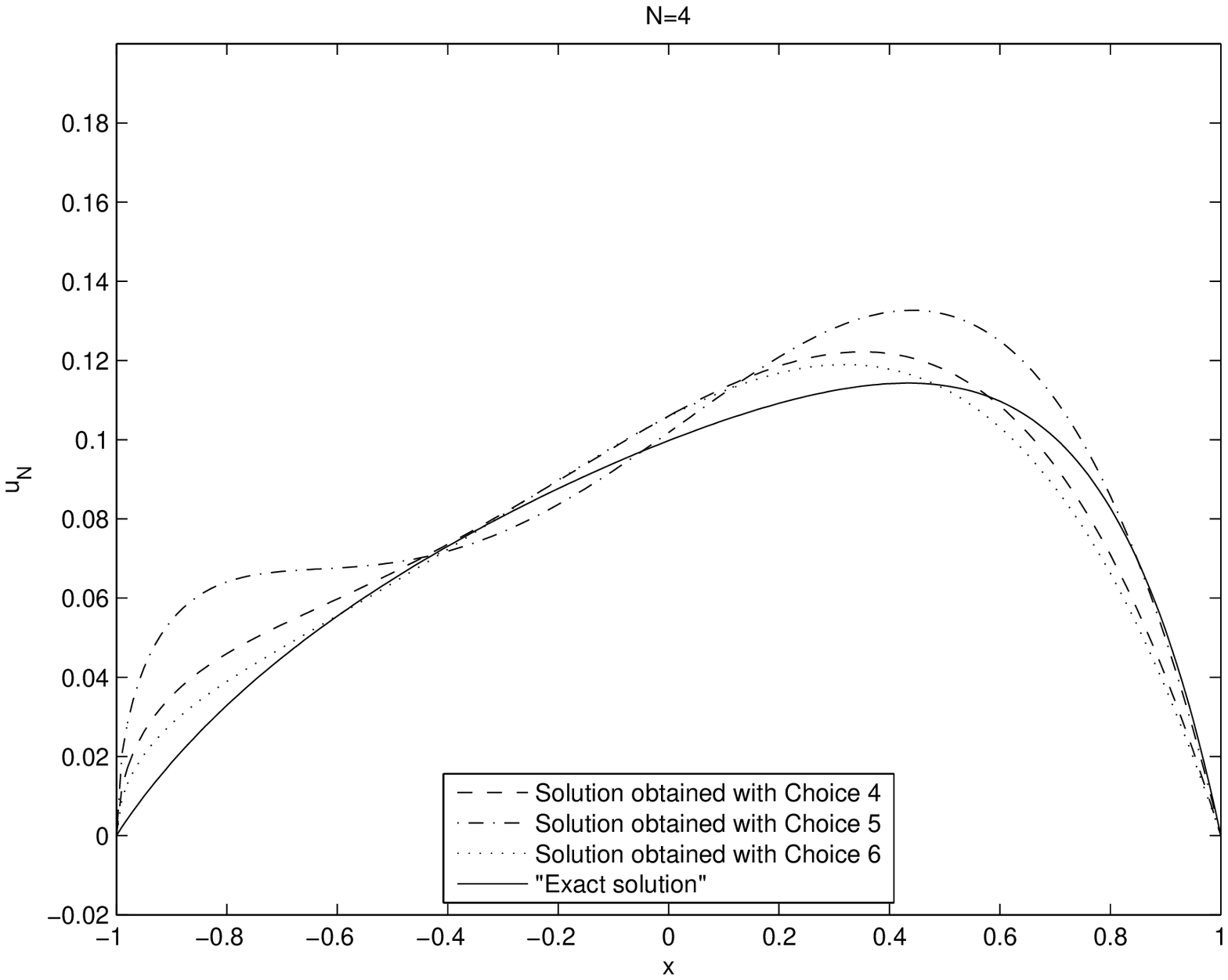}
\hspace{.3cm}\includegraphics[width=5.8cm,height=4.8cm]{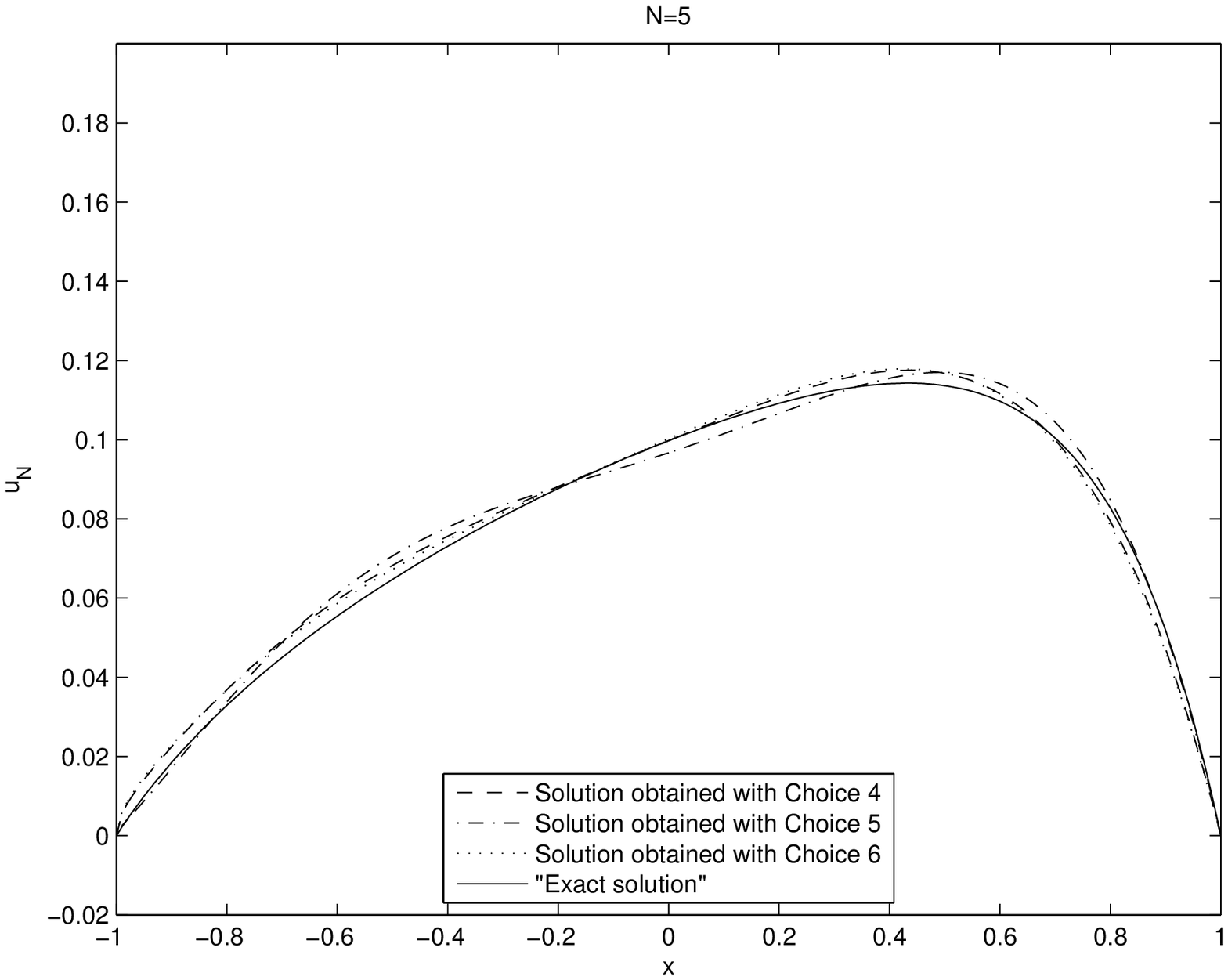}}
\caption{\small \sl 
Approximated solutions for $N=4$ and $N=5$ of the fractional differential  problem (\ref{dife}) 
with $\sigma=0.5$, $K=10$,  $g(x)=1$, for the three  different choices of  collocation 
nodes introduced in this section.}  \label{fig8}
\end{figure}

In the first numerical test, we discretize the fractional differential  problem
(\ref{dife}) with $\sigma=0.5$, $K=10$ and $g(x)=1$. This is a kind
of advection-diffusion problem, with a boundary layer developing on the
right-hand side. The behavior in the middle is regulated by the fractional
derivative operator. We then compare the results obtained by
implementing different sets of collocation nodes as specified above.
As done in Section 5, since the exact solution is not available, we substitute it
with an approximation obtained with $N$ sufficiently large.
Figure \ref{fig8} shows the results of this test for $N=4,5$, respectively.
The superiority of our method is evident, as  is also illustrated by the
results of Table \ref{table3} where the errors in the discrete maximum norm, relative to 
the cases examined, are shown.

\begin{table}
\centerline{
{\begin{tabular}{@{}ccclccclccclccclccc}\hline
 &&$N$ &&  Error 4&&      Error 5&&  Error 6 &&           \\
\hline
 &&    4   &&      0.0111&&  0.0276  &&     0.0045 \\
 &&     5 &&   0.0039&&      0.0049 &&  0.0040\\
&&6   && 0.0049    &&    0.0073&&   0.0030\\
&&7   && 0.0032    &&   0.0033  &&   0.0024 \\
&&8   && 0.0028    &&     0.0033  && 0.0019 \\
&&9   && 0.0022    &&     0.0025  &&  0.0016 \\
&&10 && 0.0018    &&    0.0020   &&  0.0013\\
&&11 &&0.0015     &&   0.0017  &&  0.0011\\
&&12 && 0.0013    &&  0.0014  &&   9.2276e-04\\
&&13  && 0.0011   &&  0.0012   && 7.8751e-04\\
&&14&& 9.5253e-04 &&  0.0010 &&  6.7687e-04 \\
&&15 && 8.2687e-04 &&     8.9419e-04&& 5.8520e-04\\
    \hline
  \end{tabular}}
  }\vspace{-0.5cm}
  \bigskip\caption{\small \sl Errors in 
the discrete maximum norm between the exact solution  $u$ and the approximated solution $u_{N}$
of problem (\ref{dife}) with $\sigma=0.5$, $K=10$ and $g(x)=1$. They are  obtained with
the same representation nodes, but with different collocation nodes, as a consequence of
{\sl Choice 4} (Error 4), {\sl Choice 5} (Error 5), and {\sl Choice 6} (Error 6).}
\label{table3}
\end{table}

\par\smallskip

\begin{figure}[h!]
\vspace{.1cm}
\centering
\centering{\includegraphics[width=5.8cm,height=4.8cm]{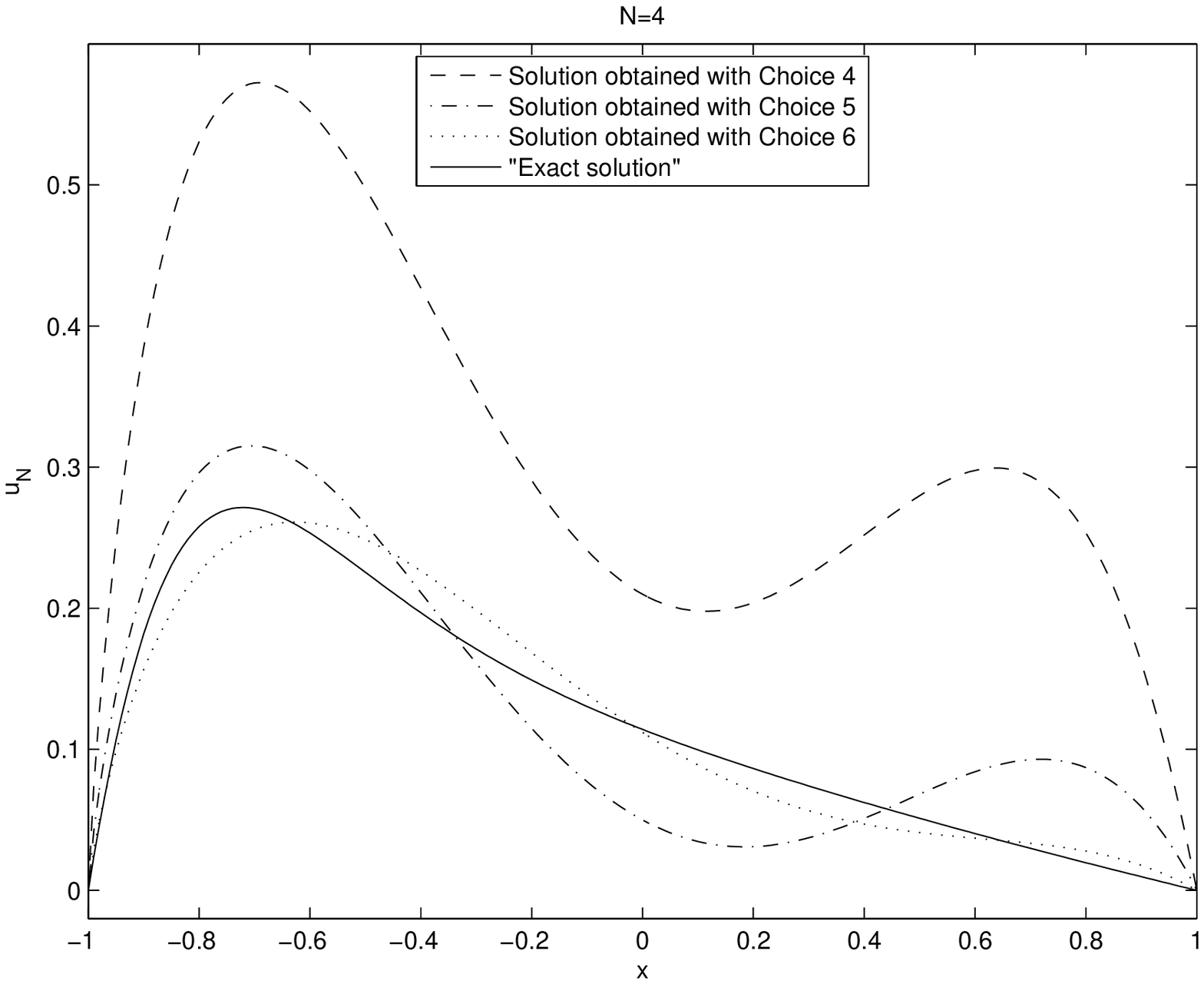}
\hspace{.3cm}\includegraphics[width=5.8cm,height=4.8cm]{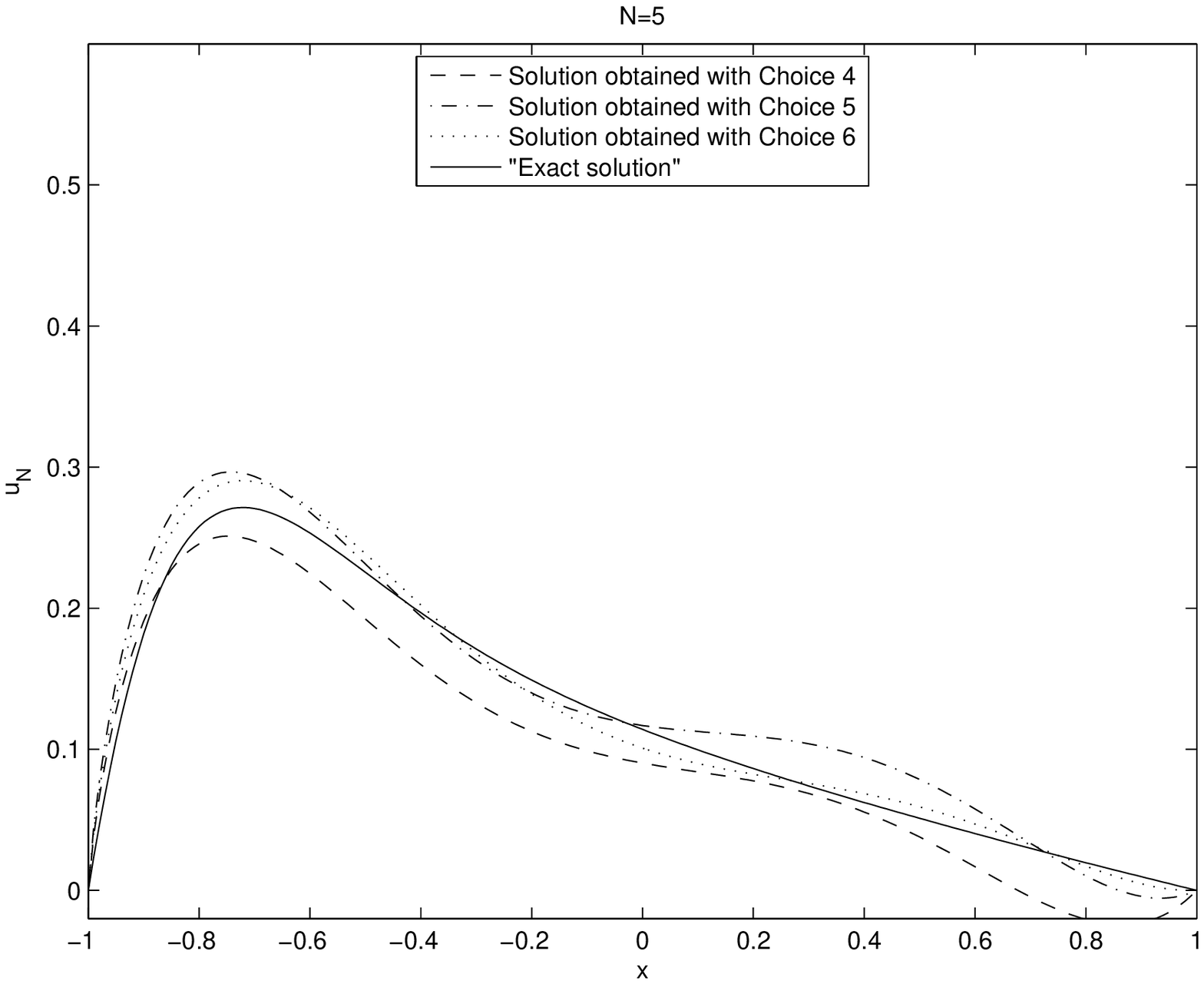}}
\caption{\small \sl 
Approximated solutions for $N=4$ and $N=5$ of the fractional differential  problem (\ref{dife}) 
with $\sigma=0.8$, $K=-10$,  $g(x)=1$, for the three  different choices of  collocation 
nodes introduced in this section.}  \label{fig11a}
\end{figure}

\begin{figure}[h!]
\vspace{.1cm}
\centering
\centering{\includegraphics[width=5.8cm,height=4.8cm]{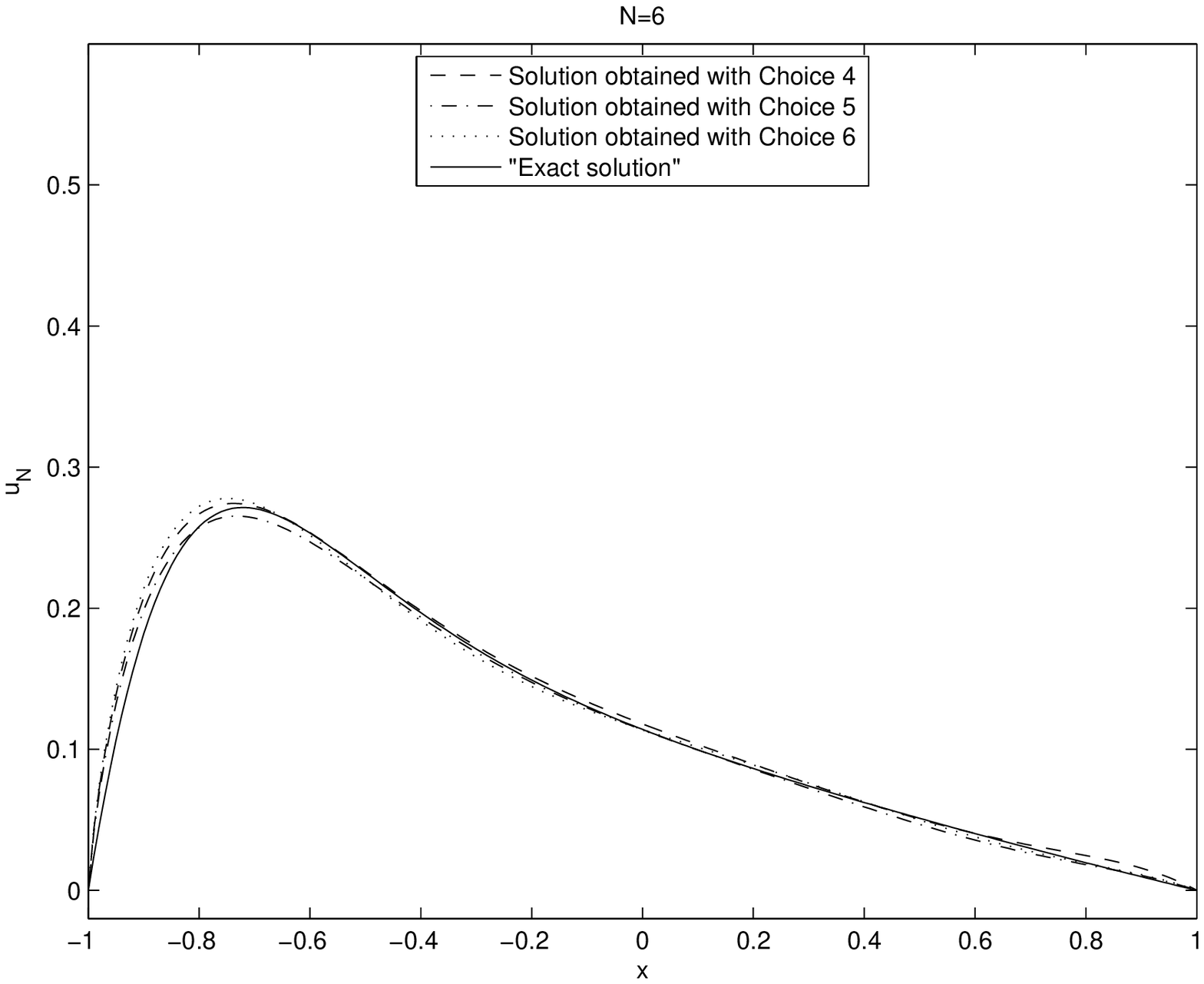}
\hspace{.3cm}\includegraphics[width=5.8cm,height=4.8cm]{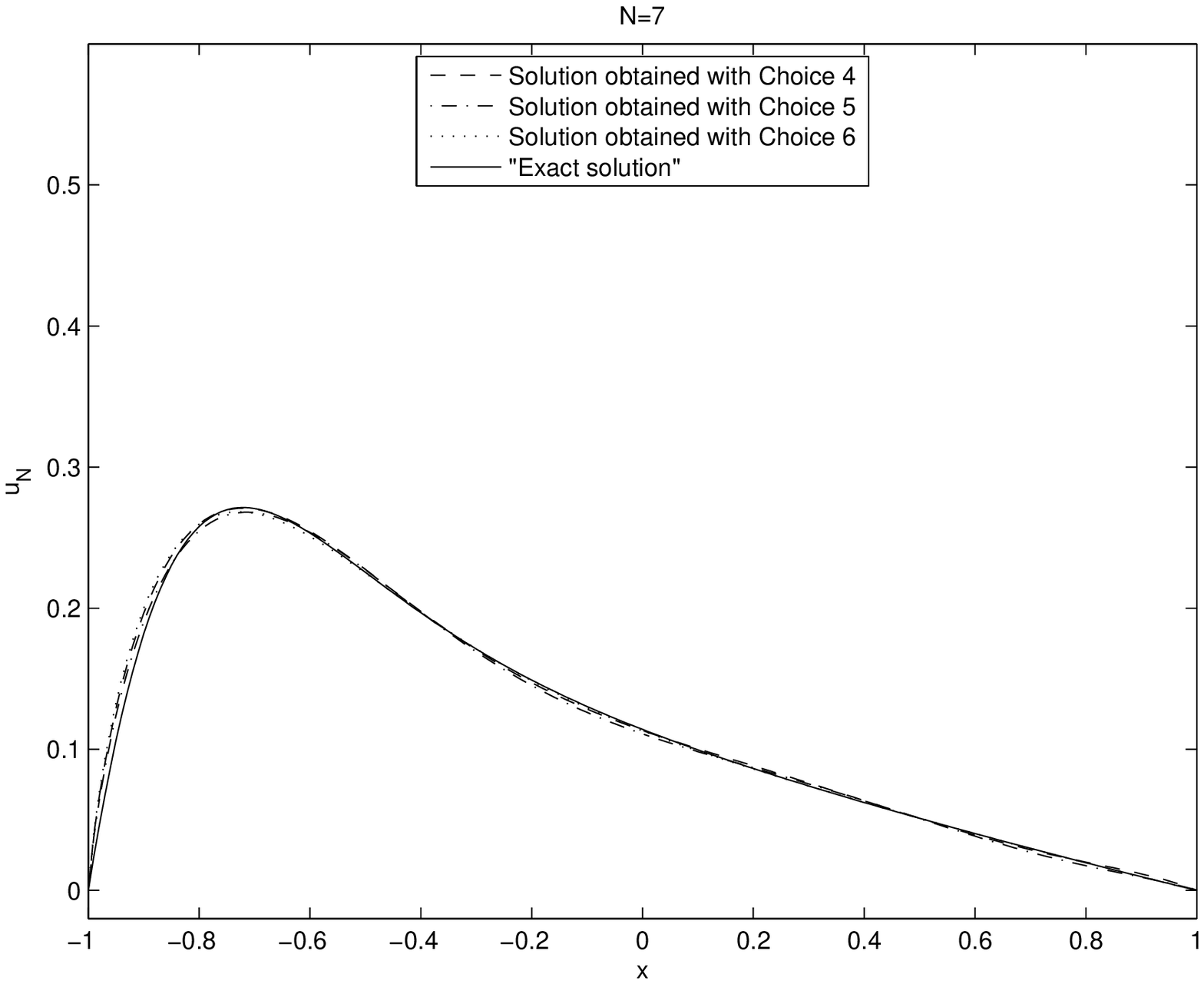}}
\caption{\small \sl 
Approximated solutions for $N=6$ and $N=7$ of the fractional differential  problem (\ref{dife}) 
with $\sigma=0.8$, $K=-10$,  $g(x)=1$, for the three  different choices of  collocation 
nodes introduced in this section.}  \label{fig11b}
\end{figure}

In the second test we have $\sigma=0.8$, $K=-10$ and $g(x)=1$.
Now, the transport is from left to right.
Figures \ref{fig11a} and \ref{fig11b} show the results for $N=4,5,6,7$.
Again, the best performances are obtained through the superconsistent
method. 
The ``exact solutions''  shown in Figures  \ref{fig8}, \ref{fig11a} and \ref{fig11b}
have been actually replaced by approximated ones, obtained for $N=50$ by
using for both representation and collocation nodes the points defined in (\ref{nodiCheb}).
\par\smallskip

Note that, for some critical values of the parameters $\sigma$, $K$, $N$,
the   superconsistent method may blow up, mainly because the procedure for computing
the collocation nodes fails. In truth, these are situations where the problem
is particularly stiff ($N$ small, $\vert K\vert$ large), a setting that 
may constitute a difficulty for any numerical technique. If one stays
within reasonable limits, our approach looks reliable and effective.
\par\smallskip

In conclusion, it may be worthwhile to spend some efforts in computing
the ``right'' set of collocation nodes (that depend on the representation
nodes and the differential operator to be approximated), since, with
no additional cost, the procedure turns out to be highly accurate and competitive,
even if compared with more standard high-order pseudo-spectral techniques.
\par

\end{document}